%% file: THHBP2Iv3.tex
\documentclass[12pt]{amsart}
\usepackage[margin=1.2in]{geometry} 
\input{THHpreamblev3}

\graphicspath{{./EU-symbol.jpg}} 
\title[THH of truncated Brown-Peterson spectra I]{Topological Hochschild homology of truncated Brown-Peterson spectra I}  

\author[Angelini-Knoll]{Gabriel Angelini-Knoll} 
\address{Universit{\'e} Sorbonne Paris Nord, LAGA, CNRS, UMR 7539, F-93430, Villetaneuse, France}
\email{angelini-knoll@math.univ-paris13.fr}   

\author[Culver]{Dominic Leon Culver}
\address{Max Planck Institute for Mathematics, Vivatsgasse 7,
53111 Bonn,   
Germany} 
\email{dculver@mpim-bonn.mpg.de}
  
\author[H\"oning]{Eva H\"oning}
\address{Radboud University, Department of Mathematics, Heyendaalseweg 135,
6525 AJ Nijmegen, The Netherlands}  
\email{eva.hoening@ru.nl}      
      
\begin{document}        
	\input{THHBP2Iabstractv3}               
    \maketitle       
	\tableofcontents     
	\input{THHBP2Iintrov3}     
 	\input{THHBP2coeffHFpv3}  
	\input{THHBP2coeffHZv3}

	\input{THHBP2coeffk1v3}

	\input{THHBP2coeffk2v3}

	\bibliographystyle{plain}
	\bibliography{THHv3}
\end{document}

%% file: THHpreamblev3.tex
\usepackage{amsmath}
\usepackage{amsthm}
\usepackage{amssymb}
\usepackage{lscape,xcolor}
\usepackage{graphicx}
\usepackage{mathrsfs}
\usepackage{mathtools}
\usepackage{stmaryrd}
\usepackage{verbatim}
\usepackage{rotating}
\usepackage{tikz-cd}
\usepackage{amsrefs}
\usepackage[colorlinks=true, linkcolor=blue, menucolor=blue]{hyperref}
\usepackage{euscript}
\usepackage[colorinlistoftodos]{todonotes}
\usepackage{spectralsequences}
\usepackage{xypic}
\usepackage{lscape}

\usepackage{xcolor}
\definecolor{seagreen}{RGB}{46,139,87}
\definecolor{maroon}{RGB}{128,0,0}
\definecolor{darkviolet}{RGB}{148,0,211}
\definecolor{twelve}{RGB}{100,100,170}
\definecolor{thirteen}{RGB}{100,150,50}
\definecolor{fourteen}{RGB}{200,0,0}
\definecolor{fifteen}{RGB}{0,200,0}
\definecolor{sixteen}{RGB}{0,0,200}
\definecolor{seventeen}{RGB}{200,0,200}
\definecolor{eighteen}{RGB}{0,200,200}

\allowdisplaybreaks[1]

\newcommand{\mmod}{\! \sslash \!}

\newcommand{\mr}[1]{\mathrm{#1}}

\newcommand{\ms}[1]{\mathscr{#1}}

\newcommand{\bra}[1]{\langle #1 \rangle}
\newcommand{\br}[1]{\overline{#1}}

\newcommand{\Z}{\mathbb{Z}}

\newcommand{\Q}{\mathbb{Q}}

\newcommand{\F}{\mathbb{F}}

\newcommand{\tBP}[1]{BP\bra{#1}}
\newcommand{\B}{\tBP{2}}
\newcommand{\Bn}{B \langle n \rangle}
\newcommand{\Bt}{B \langle 2 \rangle}

\newcommand{\tmf}{\mathrm{tmf}}
\newcommand{\taf}{\mathrm{taf}^{D}}

\def \HF2{\mr{H}\F_2}

\newcommand{\MU}{\mr{MU}}

\newcommand{\BP}{\mr{BP}}

\DeclareMathOperator{\Ext}{Ext}
\DeclareMathOperator{\Tor}{Tor}

\DeclareMathOperator*{\colim}{colim}

\DeclareMathOperator{\THH}{THH}
\DeclareMathOperator{\HH}{HH}

\newcommand{\A}{\ms{A}}

\def \AA0{\br{A \mmod A(0)}_*}
\def \AA2{A\mmod A(2)_*}

\def \AE2{A\mmod E(2)_*}
\renewcommand{\AE}[1]{A\mmod E(#1)_*}

\def\co{\colon \thinspace}

\newcommand{\op}{\mathrm{op}}
\newcommand{\cof}{\mathrm{cof}}

\newcommand{\bxi}{\bar{\xi}}
\newcommand{\btau}{\bar{\tau}}
\newcommand{\sig}[1]{\sigma{#1}}

 \newtheorem{thm}[equation]{Theorem}
 \newtheorem{cor}[equation]{Corollary}
 \newtheorem{lem}[equation]{Lemma}
 \newtheorem{prop}[equation]{Proposition}
 
  \newtheorem{rem}[equation]{Remark}
 
 \newtheorem*{thm*}{Theorem}
 \newtheorem*{cor*}{Corollary}
 \newtheorem*{lem*}{Lemma}
 \newtheorem*{prop*}{Proposition}
  \newtheorem*{not*}{Notation}

 \theoremstyle{definition}
 \newtheorem{defn}[equation]{Definition}
 
 \newtheorem{exs}[equation]{Examples}

 \newtheorem{conjecture}[equation]{Conjecture}
\newtheorem*{defn*}{Definition}
\newtheorem*{ex*}{Example}
\newtheorem*{exs*}{Examples}
\newtheorem*{rmk*}{Remark}
\newtheorem*{claim*}{Claim}

\numberwithin{equation}{section}
\numberwithin{figure}{section}

%% file: THHBP2Iabstractv3.tex
% root file is THHBP2.tex
\begin{abstract}
We compute topological Hochschild homology of sufficiently structured forms of truncated Brown--Peterson spectra with coefficients. In particular, we compute $\THH_*(B\langle n\rangle ;H\mathbb{Z}_{(p)})$ for all $n$ and $\THH_*(B\langle 2\rangle ;M)$ for $M\in  \{ k(1),k(2)\}$ where $B\langle n\rangle$ is an $E_3$ form of $BP\langle n\rangle$ for certain primes $p$. For example, this gives a computation of $\THH(\taf;M)$ for $M\in \{H\mathbb{Z}_{(3)},k(1),k(1))$ where $\taf$ is the $E_{\infty}$ form of $BP\langle 2\rangle$ constructed by Hill--Lawson. 
\end{abstract}

%% file: THHBP2Iintrov3.tex
% root file is THHBP2.tex

\section{Introduction}
Topological Hochschild homology and cohomology are rich invariants of rings, or more generally ring spectra, with applications to such fields as string topology \cite{CohenJones}, deformation theory of $A_{\infty}$ algebras \cite{Ang08}, and integral $p$-adic Hodge theory \cite{BMSII}. Topological Hochschild homology is also a first order approximation to algebraic K-theory in a sense made precise using Goodwillie calculus by \cite{DundasMcCarthy}. 

Algebraic K-theory of ring spectra that arise in chromatic stable homotopy theory are of particular interest because of the program of Ausoni--Rognes \cite{AR02}, which suggests that algebraic K-theory shifts chromatic complexity up by one, a higher chromatic height analogue of conjectures of Lichtenbaum \cite{Lichtenbaum} and Quillen \cite{Qui75}. A higher chromatic height analogue of one of the Lichtenbaum--Quillen conjectures was recently proven for truncated Brown--Peterson spectra $\tBP{n}$ by \cite{HahnWilson18}. However, it is still desirable to have a more explicit computational understanding of algebraic K-theory of $\tBP{n}$ in order to understand the \'etale cohomology of $\tBP{n}$ as suggested by Rognes  \cite[\S 5-6]{Rog14}. 

One of the most fundamental objects in chromatic stable homotopy theory is the Brown--Peterson spectrum $\BP$, which is a complex oriented cohomology theory that carries the universal $p$-typical formal group. The coefficients of $\BP$ are the symmetric algebra over $\mathbb{Z}_{(p)}$ on generators $v_i$ for $i\ge 1$, and we may form truncated versions of $\BP$, denoted $\tBP{n}$ by coning off the regular sequence $(v_{n+1},v_{n+2}, \ldots )$. More generally, we consider forms of $\tBP{n}$, in the spirit of \cite{Morava}, which are constructed by coning off some sequence $(v_{n+1}^{\prime},v_{n+2}^{\prime},\ldots )$ of indecomposable algebra generators in $\BP_*$ where $|v_{k}^{\prime}|=|v_k|$ (see Definition \ref{form} for a precise definition). We will be most interested in working with forms of $\tBP{n}$ that are $E_m$ ring spectra for sufficiently large $m$. We will refer to such spectra as $E_m$ forms of $\tBP{n}$. For example, the spectrum $H\mathbb{Z}_{(p)}$ is an $E_{\infty}$ form of $\tBP{0}$, and $\ell$ is an $E_{\infty}$ form of $\tBP{1}$ at all primes by \cite{BakerRichter}.  
 
In the last decade $E_{\infty}$ forms of $\B$ were constructed at the prime $p=2$ by \cite{LawsonNaumann} and $p=3$ by \cite{HillLawson}. In \cite{LawsonNaumann}, Lawson--Naumann use the the moduli stack of formal groups with a $\Gamma_1(3)$-structure to construct an $E_{\infty}$ form of $\B$ at the prime $2$ denoted $\tmf_1(3)$. In \cite{HillLawson}, Hill--Lawson use a quaternion algebra $D$ of discriminant $14$ and its associated Shimura curve $\mathcal{X}^D$ to construct an $E_{\infty}$ form of $\B$ at the prime $p=3$, denoted $\taf$. Even more recently, Hahn--Wilson \cite{HahnWilson20} construct an $E_3$ form of $\tBP{n}$ at all primes and for all $n$, which we denote $\tBP{n}^{\prime}$. 
This is especially interesting since no $E_{2(p^2+2)}$ form of $\tBP{n}$ exists for $n\ge 4$ by Lawson \cite{Law18} at the prime $p=2$ and Senger \cite{Sen17} at primes $p>2$. Highly structured models for truncated Brown--Peterson spectra make computations of invariants of these truncated Brown--Peterson spectra more tractable, and therefore they will be important for our calculations. 
 
For small values of $n$, the calculations of $\THH_*(\tBP{n})$ are known and of fundamental importance. The first known computations of topological Hochschild homology are B\"okstedt's calculations of  $\THH_*(H\mathbb{F}_p)$ and $\THH_*(H\mathbb{Z}_{(p)})$ in \cite{Bok85}. To illustrate how fundamental these computations are, we point out that the computation 
$\THH_*(H\mathbb{F}_p) \cong P(\mu_0)$
where $|\mu_0|=2$ is the linchpin for a new proof of Bott periodicity \cite{HN19}. In \cite{McClureStaffeldt}, McClure and Staffeldt compute the Bockstein spectral sequence 
\[ \THH_*(\ell; H\mathbb{F}_p)[v_1]\implies \THH_*(\ell ; k(1)).\]

This result is extended in \cite{AHL} where Angeltveit, Hill, and Lawson compute the square of spectral sequences 
\[ 
\xymatrix{
\THH_*(\tBP{1}; H\mathbb{F}_p)[v_0,v_1] \ar@{=>}[r] \ar@{=>}[d] & \THH_*(\tBP{1}; H\mathbb{Z}_{(p)})_p[v_1] \ar@{=>}[d] \\
\THH_*(\tBP{1}; k(1))[v_0] \ar@{=>}[r] & \THH_*(\tBP{1};\tBP{1})_p.
}
\]
This gives a complete computation of $\THH_*(\tBP{1})$. 

Let $\Bn$ denote an $E_{3}$ form of $BP\langle n\rangle$ (see Definition \ref{form}).\footnote{Note that there is a spectrum commonly denoted $B(n)=v_n^{-1}P(n)$ in other references (e.g. \cite{greenbook}) and our notation and meaning is distinct.} We compute
\[ \THH_*(\Bn;H\F_p)\cong E(\lambda_1,\dots,\lambda_{n+1})\otimes P(\mu_{n+1}) \]  
where $|\lambda_i|=2p^i-1$ and $|\mu_{n+1}|=2p^{n+1}$ in Proposition \ref{coef fp}  as a consequence of work of \cite{AngeltveitRognes}. In \cite{HahnWilson18}, Hahn and Wilson calculate the groups $\THH_*(\Bn /\MU)$, but working over $\MU$ significantly simplifies the calculation. In \cite{AusoniRichter}, Ausoni and Richter compute $\THH_*(E(2))$ under the assumption that $E(2)=\B [v_2^{-1}]$ has an $E_{\infty}$ ring structure and give a conjectural answer for $\THH_*(E(n))$, which is consistent with our calculations. These are currently the only known results for $n\ge 2$. 

The main three results of this paper are computations of the Bockstein spectral sequences 
\begin{align}
	\label{vo bockstein} \THH_*(\Bn ;H\mathbb{F}_p)[v_0]& \implies \THH_*(\Bn ; H\mathbb{Z}_{(p)})_p\\ 
	\label{v1 bockstein} \THH_*(B\langle 2\rangle ;H\mathbb{F}_p)[v_1]&\implies \THH_*(B\langle 2\rangle ; k(1))\\ 
	\label{v2 bockstein}  \THH_*(B\langle 2\rangle ;H\mathbb{F}_p)[v_2]&\implies \THH_*(B\langle 2\rangle ; k(2))
\end{align}
where $\Bn$ is an $E_3$ form of $BP\langle n\rangle$ and we assume $p\ge 3$ for our computation of the spectral sequence \eqref{v1 bockstein}. The Bockstein spectral sequences \eqref{v1 bockstein} and \eqref{v2 bockstein} are of similar computational complexity to the main result of McClure--Staffeldt \cite{McClureStaffeldt} and we were inspired by their work. 

We summarize our three main results as follows: First, we compute the topological Hochschild homology of an $E_3$ form of $\tBP{n}$ with $H\mathbb{Z}_{(p)}$ coefficients. 
\begin{thm*}[Theorem \ref{BHZ}]
Let $\Bn$ be an $E_3$ form of $\tBP{n}$ and at $p>2$ assume the error term \eqref{eq: obstruction} 
vanishes. Then there is an isomorphism of graded $\Z_{(p)}$-modules
\[ \THH_*(\Bn ;H\Z_{(p)})\cong E_{\Z_{(p)}}(\lambda_1,\ldots,\lambda_{n}) \otimes  \left (\Z_{(p)}\oplus T_0^n\right )\]
where $T_0^n$ is an explicit torsion $\Z_{(p)}$-module defined in \eqref{eq: t0n}.
\end{thm*}
In particular, the error term \eqref{eq: obstruction} vanishes for any $E_4$ form of $BP\langle n\rangle$ such as $B\langle 2\rangle=\taf$.  It is possible that the error term \eqref{eq: obstruction} also vanishes for $B\langle n\rangle  =\tBP{n}^{\prime}$ where $\tBP{n}^{\prime}$ is the $E_3$ form of $\tBP{n}$ constructed by Hahn--Wilson \cite{HahnWilson20} at odd primes, but it is not known to the authors. Theorem \ref{BHZ} also holds for $B\langle 2\rangle=\tmf_1(3)$ and $\Bn =\tBP{n}^{\prime}$, where $\tBP{n}^{\prime}$ is the $E_3$ form of $\tBP{n}$ at the prime $2$ constructed by Hahn--Wilson \cite{HahnWilson20}.  

Second, we compute the topological Hochschild homology  of an $E_{3}$ form $B\langle 2\rangle$ of $\B$ at $p\ge 3$ with $k(1)$ coefficients. 
\begin{thm*}[Theorem \ref{mod p v_2}]
Let $B\langle 2\rangle$ denote an $E_{3}$ form of $\B$ at an odd prime $p$. There is an isomorphism of $P(v_1)$-modules 
\[ \THH_*(B\langle 2\rangle;k(1))\cong E(\lambda_1)  \otimes  \left ( P(v_1)\oplus T_1^2 \right )\]
where $T_1^2$ is an explicit $v_1$-torsion $P(v_1)$-module defined in \eqref{eq: t12}. 
\end{thm*}
In particular, this result holds for $B\langle 2\rangle=\taf$ and $BP\langle 2\rangle^{\prime}$ at odd primes. 

Finally, we compute topological Hochschild homology of any $E_{3}$ form of  $\B$ with $k(2)$ coefficients.
\begin{thm*}[Theorem \ref{mod p v_1}]
Let $B\langle 2\rangle$ be an $E_{3}$ form of $BP\langle 2\rangle$. 
There is an isomorphism of $P(v_2)$-modules 
\[ \THH_*(B\langle 2\rangle;k(2))\cong P(v_2) \oplus T_2^2\]
where $T_2^2$ is an explicit $v_2$-torsion $P(v_2)$-module defined in \eqref{eq: t22}. 
\end{thm*}
\noindent 
In particular, this result holds for $B\langle 2\rangle =\taf$,  $B\langle 2\rangle=\tmf_1(3)$, and $BP\langle 2\rangle^{\prime}$ at any prime. We end with a conjectural answer (cf. Conjecture \ref{conj}) for $\THH_*(\Bn ; k(m))$ for all integers $1\le m\le n$ and any $E_3$ form of $\Bn$ at a prime $p$. 

We now outline our approach to computing $\THH_*(\taf)$ in the sequels to this paper. There is a cube of Bockstein spectral sequences
\begin{align}\label{bockstein cube}
\xymatrix{ 
 H\mathbb{F}_3 \ar@{=>}[rd] \ar@{=>}[dd] \ar@{=>}[rr] & & H\mathbb{Z}_{(3)} \ar@{=>}[rd] \ar@{=>}[dd] & \\
 & k(1) \ar@{=>}[dd] \ar@{=>}[rr] & & B\langle 2\rangle /v_2\ar@{=>}[dd]  \\
k(2) \ar@{=>}[dr] \ar@{=>}[rr] & & B\langle 2\rangle/v_1 \ar@{=>}[dr] & \\
 &  B\langle 2\rangle  /3 \ar@{=>}[rr] & &  B\langle 2\rangle  
 }   
\end{align}
where we use the abbreviation $M/x\implies M$ for the Bockstein spectral sequence with signature
\[ \THH_*(\taf ;M/x)[x]\implies \THH_*(\taf ;M)\]
where 
$M\in \{H\mathbb{Z}_{(3)},k(1),k(2),\taf/3,\taf/v_1,\taf/v_2,\taf\}.$
Here we write $\taf/x$ for the cofiber of a representative of an element $x\in \pi_{2k}\taf$ regarded as a $\taf$-module map $\Sigma^{2k}\taf\to \taf$. In the sequels to this paper, we plan to compute $\THH_*(\taf ;M)$ for $M=\taf/3$ and $M=\taf/v_1$ by comparing the edges of the cube of Bockstein spectral sequences to the Hochschild--May spectral sequence \cite{AKS18} and the Brun spectral sequence \cite{Hon20}, which compute the diagonals of the faces of the cube directly. Finally, we plan to compute $\THH_*(\taf)$ by again comparing the Hochschild--May spectral sequence to the relevant Bockstein spectral sequences in addition to cosimplicial descent techniques. 
 
\subsubsection*{Conventions}
We write $F_*X$ for $\pi_*(F\wedge X)$ for any spectra $F$ and $X$.
We also use the shorthand $H_*(X)$ for $(H\mathbb{F}_p)_*X$ for any spectrum $X$.   
 We write $\dot{=}$ to mean that an equality holds up to multiplication by a unit. The dual Steenrod algebra $H_*(H\mathbb{F}_p)$ will be denoted $\A_*$ with coproduct $\Delta\co \A_*\to \A_*\otimes \A_*$. Given a left $\A_*$-comodule $M$, its left coaction will be denoted $\nu\co \A_*\to \A_*\otimes M$ where the comodule $M$ is understood from the context. The antipode $\chi\co \A_*\to\A_*$ will not play a role except that we will write $\bar{\xi}_i:=\chi(\xi_i)$ and $\bar{\tau}_i:=\chi(\tau_i)$. 

When not otherwise specified, tensor products will be taken over $\mathbb{F}_p$ and $\HH_*(A)$ denotes the Hochschild homology of a graded $\mathbb{F}_p$-algebra relative to $\mathbb{F}_p$. We will let $P_R(x)$, $E_R(x)$ and $\Gamma_R(x)$ denote a polynomial algebra, exterior algebra, and divided power algebra over $R$ on a generator $x$. When $R=\mathbb{F}_p$, we omit it from the notation. Let $P_i(x)$ denote the truncated polynomial algebra $P(x)/(x^i)$. 

\subsubsection*{Acknowledgements}
The authors would like to thank Tyler Lawson for helpful conversations. The authors would also like to thank John Rognes and an anonymous referee for their careful reading of the paper and useful suggestions for improvement. 
Parts of this paper were written while the second author was in residence at the Max Planck Institute for Mathematics. He would like to thank the Institute for their support. The third author  was funded by the DFG priority
program SPP 1786 Homotopy Theory and Algebraic Geometry and the Radboud Excellence Initiative. This project has received funding from the European Union's Horizon 2020 research and innovation programme under the Marie Sk\l{}odowska-Curie grant agreement No 1010342555.
\thinspace \includegraphics[scale=0.1]{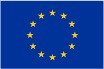} 

%% file: THHBP2coeffHFpv3.tex
% root file is THHBP2BP1.tex

\section{Topological Hochschild homology mod $(p,\dots,v_n)$}\label{sec prelim}
We begin by giving a precise definition of an $E_m$ form $\Bn$ of $\tBP{n}$. We then compute topological Hochschild homology of an $E_3$ form $\Bn$ of $BP\langle n\rangle$ at an arbitrary prime $p$ with coefficients in $H\mathbb{F}_p$. First, recall that there is an isomorphism $BP_*\cong \mathbb{Z}_{(p)}[v_i:i\ge 1]$ and an isomorphism
\[\BP_*\BP\cong \mathbb{Z}_{(p)}[v_i : i\ge 1][t_i : i\ge 1]\]
where the degrees of the generators are $|v_i|=|t_i|=2p^i-2$ for $i\ge 1$. The generators $t_i$ are determined by the canonical strict isomorphism $f$ from the universal $p$-typical formal group law to itself given by the power series 
\[ f^{-1}(x)=\sum_{i\ge 0}^{F} t_ix^{p^i}\]
where $F$ is the universal $p$-typical formal group law \cite[Lemma A2.1.26]{greenbook}. We let $v_i$ be the Araki generators. Note 
that the Araki generators agree with Hazewinkel generators mod $p$ \cite[Theorem A2.2.3]{greenbook}.  
\input{Formsv3}

\subsection{Topological Hochschild homology mod $(p,\dots ,v_n)$}
The mod $p$ homology of $\THH(\tBP{n})$ has been calculated by Angeltveit--Rognes in \cite[Theorem 5.12]{AngeltveitRognes} assuming that $\tBP{n}$ is an $E_3$ ring spectrum.  Their argument also applies to topological Hochschild homology of any $E_3$ form $B\langle n\rangle$ of $\tBP{n}$ at a prime $p$, as we now explain.  
By Proposition \ref{properties of forms}, the linearization map \eqref{map to HFp} induces an isomorphism
\[H_*(B\langle n \rangle)\cong \begin{cases} P(\bxi_1,\bxi_2,\ldots )\otimes E(\btau_{n+1},\btau_{n+2},\ldots)  & \text{ if } p \ge   3 \\ P(\bxi_1^2,\dots \bxi_{n+1}^2,\bxi_{n+2},\ldots) & \text{ if } p=2 \end{cases}\]
with its image in $\A_*$ as sub $\A_*$-comodule algebra of $\A_*$.  
By \cite[Theorem 3.4]{BFV}, the spectrum $\THH(\Bn; H\mathbb{F}_p)$ is an $E_2$ ring spectrum and the unit map 
\[H\mathbb{F}_p \to \THH(\Bn ; H\mathbb{F}_p)\] 
is a map of $E_2$ ring spectra. Using \cite[\S 3.3]{BFV} the proof of \cite[Proposition 4.3]{AngeltveitRognes} carries over mutatis mutandis and implies that the B\"okstedt spectral sequence with signature
\[ E^2_{*,*} = \HH_{*,*}(H_*(\Bn); \A_*) \Longrightarrow H_*(\THH(\Bn; H\mathbb{F}_p))\] 
is a spectral sequence of $\A_*$-comodule algebras. As in \cite[\S 5.2]{AngeltveitRognes}, the spectral sequence collapses at the $E^2$-page if $p = 2$. 
If $p \geq 3$, one can use the map to the B\"okstedt spectral sequence with signature 
\[ E^2_{*,*} = \HH_{*,*}(\A_*) \Longrightarrow H_*(\THH(H\mathbb{F}_p))\]
to determine the differentials (cf. \cite[\S 5.4]{AngeltveitRognes}). Since $\Bn$ is an $E_3$ ring spectrum, Dyer--Lashof operations are defined on $H_*(\Bn)$ and $H_*(\THH(\Bn; H\mathbb{F}_p))$ in a range that is sufficient to resolve the multiplicative extensions (see \cite[Proof of Theorem 5.12]{AngeltveitRognes}). We get  an isomorphism of $\A_*$-comodule $\A_*$-algebras
\begin{align}\label{HTHHB} 
H_*(\THH(B\langle n \rangle; H\mathbb{F}_p))\cong 
	\begin{cases} 
	\A_*\otimes E(\sig{\bxi_1},\dots,\sig{\bxi_{n+1}})\otimes P(\sig{\btau_{n+1}}) &\text{ if } p \ge 3 \\ 
	\A_* \otimes E(\sig{\bxi^2_1},\dots,\sig{\bxi^2_{n+1}})\otimes P(\sig{\bxi_{n+2}}) & \text{ if }p=2. 
	\end{cases}
\end{align}
 Since $\sigma\colon \thinspace H_*(\Bn)  \to H_{*+1}(\THH(\Bn)) \to H_{*+1}(\THH(\Bn; H\mathbb{F}_p))$ is a comodule map and a  derivation,  the $\A_*$-coaction of 
\[H_*(\THH(\Bn;H\mathbb{F}_p))\]  
can be deduced from that of $H_*(\Bn ) \subseteq \A_*$ (cf. \cite[Proof of Theorem 5.12]{AngeltveitRognes}): for $p \ge 3$ the classes $\sig \bxi_i$ for $1 \leq i \leq n+1$ are $\A_*$-comodule primitives and  we have 
\begin{equation}\label{coactp>2}
\nu ( \sig{\btau_{n+1}})= 1\otimes \sig{\btau}_{n+1}+\btau_0\otimes \sig{\bxi_{n+1}}. 
\end{equation}
For $p = 2$ the classes $\sig \bxi^2_i$  for $1 \leq i \leq n+1$  are $\A_*$-comodule primitives and we have 
\begin{equation}\label{coactp=2} \nu(\sig{\bxi_{n+2}})=1\otimes \sig{\bxi_{n+2}}+ \bxi_1 \otimes \sig{\bxi_{n+1}^2}.
\end{equation}
\begin{prop}\label{coef fp}
Let $\Bn$ be an $E_3$ form of $\tBP{n}$. 
There is an isomorphism of graded $\mathbb{F}_p$-algebras
\begin{equation}\label{eqn:THH(;F_p)}
\THH_*(\Bn;H\F_p)\cong E(\lambda_1, \dots, \lambda_{n+1})\otimes P(\mu_{n+1}), 
\end{equation}
where the degrees of the algebra generators are $|\lambda_i|=2p^i-1$ for $1\le i\le n+1$ and $|\mu_{n+1}| = 2p^{n+1}$. 
\end{prop}
\begin{proof}
Since $\THH(\Bn;H\F_p)$ is an $H\F_p$-module,  the Hurewicz homomorphism  induces an isomorphism between $\THH_*(\Bn;H\F_p)$ and the  subalgebra of comodule primitives in 
$H_*(\THH(\Bn;H\F_p))$. 
For $1 \le i \le n+1$ we write $\lambda_i \coloneqq \sigma\bxi_i$ if  $p\ge 3$ and $\lambda_i \coloneqq \sigma \bxi_i^2$ if $p=2$. We also define 
\[ \mu_{n+1} \coloneqq \begin{cases} \sigma \btau_{n+1} -\btau_0 \sigma \bxi_{n+1} &  \text{ if }  p\ge 3 \\ \sigma \bxi_{n+2}-\bxi_1 \sigma \bxi^2_{n+1} & \text{ if } p=2. \end{cases} \]
Then it is clear that the subalgebra of 
$H_*(\THH(\Bn ;H\mathbb{F}_p))$ consisting of comodule primitives
is as claimed. 
\end{proof}

%% file: Formsv3.tex
% root file is THHBP2BP1.tex

\subsection{Forms of $\tBP{n}$}
We fix a precise notion of a form of the truncated Brown--Peterson spectrum in the spirit of \cite{Morava} below. 

\begin{defn}[cf. {\cite[Definition 4.1]{LawNau}}]\label{form}
Fix integers $m\ge 1$ and $n\ge 0$. By an $E_m$ form of $\tBP{n}$ (at the prime $p$), we mean a $p$-local  $E_m$ ring spectrum $R$ equipped with a complex orientation $\MU_{(p)} \to R$ such that the composite
\[ \mathbb{Z}_{(p)}[v_1,\dots ,v_n]\to \BP_*\to \pi_*\MU_{(p)}\to \pi_*R \]
is an isomorphism.  
\end{defn}
\begin{rem}
Note that we do not assume that an $E_m$ form of  $\tBP{n}$ at the prime $p$ is an $E_m$-$MU$-algebra and therefore Definition \ref{form} differs slightly from the definition of an $E_m$-$\MU$-algebra form of $\tBP{n}$ appearing in work of Hahn--Wilson \cite[Definition 2.0.1]{HahnWilson20}. A $E_m$-$\MU$-algebra form of $\tBP{n}$  in the sense of \cite[Definition 2.0.1]{HahnWilson20} is a $E_m$ form of $\tBP{n}$ at the prime $p$ in the sense of Definition \ref{form}. The distinction arises because for example $\taf$ is an $E_{\infty}$ form of $BP\langle 2\rangle$, however it is not known, at least to the authors, whether the complex orientation $\MU\to \taf$ can be elevated to an $E_{\infty}$-ring spectrum map. Nonetheless, we know that the map $\MU\to \taf$ is an $E_2$-ring spectrum map by \cite[Theorem 1.2]{CM15}, which is sufficient for our purposes. 
\end{rem}

\begin{not*}
Throughout, we let $\Bn$ denote an $E_3$ form of $\tBP{n}$ at the prime $p$ in the sense of Definition \ref{form} for $n\ge 0$. 
\end{not*}
We collect some consequences of Definition \ref{form}.

\begin{prop}\label{properties of forms}
Since $\Bn$ is an $E_3$ form of $\tBP{n}$ at the prime $p$ for $m\ge 3$, the following hold: 
\begin{enumerate}
\item \label{part 1 form lemma} There are indecomposable algebra generators $v_i^{\prime}$ with $v_i^{\prime}=v_i$ for $1\le i\le n$ such that 
$\BP_*/(v_k^{\prime} : k\ge n+1)\cong \pi_*B \langle n \rangle$.
\item \label{part 2 form lemma}The orientation 
$\MU_{(p)}\to \Bn$ lifts to an $E_2$ ring spectrum map and consequently there is an $E_2$ ring spectrum map $\BP\to \Bn$ realizing the canonical quotient map $\BP_*\longrightarrow \BP_*/(v_k^{\prime} : k\ge n+1)$ on homotopy groups.   
\item \label{part 3 form lemma} There is an $E_3$ ring spectrum map $\Bn\to H\mathbb{Z}_{(p)}$ and the map 
induced by the composite 
\begin{align}\label{map to HFp}
	\Bn \to H\mathbb{Z}_{(p)}\to H\mathbb{F}_p
\end{align}
in mod $p$ homology provides an isomorphism 
\[ H_*(\Bn )\cong \A//E(n)_*\subset \A_*\]
of $\A_*$-comodule $\mathbb{F}_p$-algebras onto its image in the dual Steenrod algebra. 
\item \label{part 4 form lemma}  If $B \langle n \rangle$ is $E_{3}$ and $x_1, \dots , x_n$ is a regular sequence of elements in $B \langle n \rangle_*$, then one can construct the spectrum 
 $\Bn/(x_1,x_2,\dots ,x_n)$ as an $E_1$ $B \langle n \rangle$-algebra.
\item \label{part 5 form lemma} The $p$-completion of $\Bn$ is weakly equivalent to the $p$-completion of any other $E_m$ form of $\tBP{n}$ at the prime $p$ in the category of spectra.
\end{enumerate}
\end{prop}
\begin{proof}
For Part \eqref{part 1 form lemma} set $v_i' \coloneqq v_i - f_i(v_1, \dots, v_n)$ for $i \geq n+1$, where $f_i(v_1, \dots, v_n)$ is the image of $v_i$ under $\BP_* \to \tBP{n}_*  \cong \mathbb{Z}_{(p)}[v_1, \dots, v_n]$. 
Part \eqref{part 2 form lemma} follows by applying \cite[Theorem 1.2]{CM15}. 
Part \eqref{part 3 form lemma} is \cite[Theorem 4.4]{LawNau}. 
Part \eqref{part 4 form lemma} follows from \cite[Section 3]{Ang08} (cf. \cite[Theorem A]{HahnWilson18}).
Part \eqref{part 5 form lemma} is \cite[Theorem A]{AL17}. 
\end{proof}

\begin{exs}
 The Eilenberg--MacLane spectrum $H\mathbb{Z}_{(p)}$ is an $E_{\infty}$ form of $\tBP{0}$. The Adams summand $\ell$ is an $E_{\infty}$ form of $\tBP{1}$ by \cite[Corollary 1.4]{BakerRichter}. 
\end{exs}

\begin{not*}
Let $\tmf_1(3)$ denote the $E_{\infty}$ form of $\B$ constructed by Lawson--Naumann \cite{LawsonNaumann} at $p=2$. Let $\taf$ denote the $E_\infty$ form of $\B$ constructed by Hill--Lawson \cite{HillLawson} at $p=3$.	Let $\tBP{n}^{\prime}$ denote the $E_3$ form of $\tBP{n}$ constructed by Hahn--Wilson \cite{HahnWilson20} at all primes. 
\end{not*}

%% file: THHBP2coeffHZv3.tex
% root file is THHBP2BP1.tex

\section{Topological Hochschild homology mod $(v_1, \dots ,v_n)$}
We begin by setting up the Bockstein spectral sequence. In order to ensure that this spectral sequence is multiplicative, we compare it with the Adams spectral sequence.

\subsection{Bockstein and Adams spectral sequences}\label{Bockstein and Adams}
Let $\Bn$ be an $E_3$ form of $BP\langle n\rangle$ at the prime $p$ which is equipped with a choice of generators $v_i$ in degrees $|v_i|=2p^i-2$ for $0<i\le n$ such that $\Bn_*=\mathbb{Z}_{(p)}[v_1,\dots ,v_n]$. Let $v_0=p$ by convention. Let 
\[k(i)=\Bn/(p,\dots ,v_{i-1},v_{i+1},\dots  v_n)\] 
be the $E_1$ $\Bn$-algebra constructed in Proposition \ref{properties of forms} \eqref{part 4 form lemma} where $k(0)=H\mathbb{Z}_{(p)}$. We regard $k(i)$ as a right $B\langle n\rangle \wedge B\langle n\rangle^{\op}$-module by restriction along the map
\[B\langle n\rangle\wedge B\langle n\rangle^{\op}\to B\langle n\rangle \to k(i).\] 
For $0\le i\le n$ we have cofiber sequences of right $\Bn\wedge \Bn^{\op}$-modules 
\[
	\begin{tikzcd}
		\Sigma^{|v_i|}k(i)\arrow[r,"\cdot v_i"] & k(i) \arrow[r] & H\F_p. 
	\end{tikzcd} 
	\] 
Applying the functor $- \wedge_{\Bn\wedge \Bn^{\op}}\Bn$ produces the cofiber sequence 
	\[
	\begin{tikzcd}
		\Sigma^{|v_i|}\THH(\Bn;k(i))\arrow[r] & \THH(\Bn;k(i))\arrow[r] & \THH(\Bn;H\F_p).
	\end{tikzcd}
\]

Iterating this, we produce the following tower
\begin{equation}\label{eq: tower defining Bockstein}
	\begin{tikzcd}[column sep = small]
		\cdots \arrow[r] & \Sigma^{2|v_i|}T(k(i)) \arrow[r,"\cdot v_i"] \arrow[d]& \Sigma^{|v_i|}T(k(i)) \arrow[d]\arrow[r,"\cdot v_i"] & T(k(i)) \arrow[d]\\
						& \Sigma^{2|v_i|}T(H\F_p) & \Sigma^{|v_i|}T(H\F_p) & T(H\F_p),
	\end{tikzcd}
\end{equation}
where $T(k(i)) \coloneqq \THH(\Bn; k(i))$ and $T(H\mathbb{F}_p) \coloneqq \THH(\Bn; H\mathbb{F}_p)$. 

This yields an exact couple after applying homotopy groups and it produces the $v_i$-Bockstein spectral  sequence with $E_1$-page 
\begin{equation}\label{Bockstein vi}
E_1^{*,*} = \THH(\Bn; H\mathbb{F}_p)[v_i]
\end{equation}
Note that the fact that $\Bn$ and $k(i)$ are connective and have homotopy groups that are degreewise finitely generated $\mathbb{Z}_{(p)}$-modules implies that the homotopy groups of $\THH(\Bn; k(i))$ are degreewise finitely generated $\mathbb{Z}_{(p)}$-modules, too. 
It follows that $\THH_*(\Bn; k(i))$ has the form 
 \[ \THH_*(\Bn; k(i)) \cong \bigoplus_l P(v_i)\{\alpha_l\} \oplus \bigoplus_k P_{r_k}(v_i) \{\beta_k\}\]
 for some classes $\alpha_l$ and $\beta_k$. 
Here, for $i = 0$, $P(v_i)$ is defined to be $\mathbb{Z}_{(p)}$ and $P_r(v_i)$ is $\mathbb{Z}/{p^i}$. 
 We get that 
 \[ \THH_*(\Bn;H\mathbb{F}_p) \cong  \bigoplus_l \mathbb{F}_p \{a_l\} \oplus \bigoplus_k \mathbb{F}_p \{b_k\} \oplus \bigoplus_k\mathbb{F}_p \{c_k\}, \] 
 where $a_l$ and $b_k$ are the images of $\alpha_l$ and $\beta_k$ under the map $\THH_*(\Bn; k(i)) \to \THH_*(\Bn; H\mathbb{F}_p)$, and $c_k$ is a preimage of $v_i^{r_k-1}\beta_k$ under the map
 $\THH_*(\Bn; H\mathbb{F}_p) \to \Sigma^{|v_i|+1}THH_*(\Bt;k(i))$. 
 The differentials in the spectral sequence are given as follows: 
 The classes $a_l$ and $b_k$ are infinite cycles. The class $c_k$ survives to the $E_{r_k}$-page and we have 
 \[ d_{r_k}(c_k)  = v_i^{r_k} b_k.\] 
 The spectral sequence converges strongly to $\THH_*(\Bn; k(i))$ for 
 $0 < i \leq n$ and to $\pi_*(\THH(\Bn; H\mathbb{Z}_{(p)})_p)$ for $i = 0$. The cofibers in the tower \eqref{eq: tower defining Bockstein} are $H\F_p$-module spectra.  
 
 We now relate the Bockstein spectral sequence to the Adams spectral sequence. In order to do this, we show that the tower \eqref{eq: tower defining Bockstein} is also an Adams resolution. For the definition of an Adams resolution, the reader is referred to \cite[Definition 2.1.3]{greenbook}.  
In order to show that this tower is an Adams resolution, it must be shown that the vertical morphisms
\begin{equation}\label{eq: vertical morphism}
	\Sigma^{m|v_i|}\THH(\Bn; k(i))\longrightarrow \Sigma^{m|v_i|}\THH(\Bn; H\F_p)
\end{equation}
induce monomorphisms in mod $p$ homology. We have equivalences of spectra
\[
	\THH(\Bn; M)\simeq M\wedge_{\Bn}\THH(\Bn)
\]
for $M\in \{H\mathbb{F}_p,k(i): 0\le i\le n\}$ by \cite[Remark 6.1.4]{HahnWilson20} and consequently there is is an Eilenberg--Moore spectral sequence
\[
	\Tor^{H_*\Bn}_{*,*}(H_*(M), H_*(\THH(\Bn)))\implies H_*(\THH(\Bn;M))
\]
for each 
$M\in \{H\mathbb{F}_p,k(i):0\le i\le n\}$. Since $H_*(\THH(\Bn))$ is a free $H_*(\Bn)$-module by \cite[Theorem 5.12]{AngeltveitRognes}, the Eilenberg--Moore spectral sequence collapses at the $E^2$-page without room for differentials. Furthermore, the morphism \eqref{eq: vertical morphism} induces a morphism of Eilenberg--Moore spectral sequences. Thus, we observe that the morphism \eqref{eq: vertical morphism} induces the map 
\begin{align}\label{induced map in homology}
	H_*(k(i))\otimes_{H_*(\Bn)}H_*(\THH(\Bn))\to \A_*\otimes_{H_*(\Bn)}H_*(\THH(\Bn)),
\end{align}
in mod $p$ homology where the map on the first factor is induced by the linearization map $k(i)\to H\F_p$. The map \eqref{induced map in homology} is an injection. Since $H_*(\THH(\Bn))$ is a free $H_*(\Bn)$-module, the  map \eqref{eq: vertical morphism} induces an injection on mod $p$ homology. Thus, we have shown the following proposition. 

\begin{prop}
	The tower \eqref{eq: tower defining Bockstein} is an Adams resolution. 
\end{prop}

Thus, the Adams spectral sequence for $\THH(\Bn; k(i))$ agrees with the Bockstein spectral sequence for $0\le i\le n$.  
By \cite[Theorem 2.3.3]{greenbook}, we know that the Adams spectral sequence for $\THH(\Bn; k(i))$, and consequently the Bockstein spectral sequence, is multiplicative for $0\le i\le n$ from the $E_2$-page onwards. 
To see that the Adams spectral sequence is in fact multiplicative from the $E_1$-page onwards, we prove explicitly in the case $i=0$ that the $d_1$ differential satisfies the Leibniz rule in Lemma \ref{lem: d1 mod p Bockstein}. In the case $i>0$, we can apply a change of rings isomorphism and compute explicitly that the $E_2$-page is 
\[ \Ext_{E(Q_i)_*}^{*,*}(\mathbb{F}_p,E(\lambda_1,\dots,\lambda_{n+1})\otimes P(\mu_{n+1}))=P(v_i)\otimes E(\lambda_1,\dots,\lambda_{n+1})\otimes P(\mu_{n+1}).\]
using the coactions discussed previously on $\lambda_i$ and $\mu_{n+1}$. Consequently, when $i>0$ there are no non-trivial $d_1$ differentials. Altogether, this proves the following corollary. 
\begin{cor}
	The $v_i$-Bockstein spectral sequence computing $\THH_*(\Bn; k(i))$ in the case $i\ge 1$ and $\pi_*\THH(\Bn; H\mathbb{Z}_{(p)})_p$ in the case $i=0$ is multiplicative from the $E_1$-page onwards.
\end{cor} 
\input{THHBP2coeffHQv3}

\subsection{The $v_0$-Bockstein spectral sequence}
In this section, we compute the $v_0$-Bockstein spectral sequence with signature
\begin{align}
\label{v_0Bock} E_1^{*,*}= \THH_*(\Bn;H\F_p)[v_0]\implies \THH_*(\Bn;H\mathbb{Z}_{(p)})_p
\end{align}
where $\Bn$ is an $E_3$ form of $\tBP{n}$. 
At odd primes, we must assume that a certain error term \eqref{eq: obstruction} vanishes. 
This error term vanishes for any $E_4$ form of $BP\langle n\rangle$ at odd primes, for example $\taf$.
\begin{lem}\label{lem: d1 mod p Bockstein}
There is a differential 
\begin{align*}
d_1(\mu_{n+1}) \dot{=} v_0\lambda_{n+1}
\end{align*}
in the $v_0$-Bockstein spectral sequence (\ref{v_0Bock}) and the $d_1$ differential satisfies the Leibniz rule. 
\end{lem}
\begin{proof}
We just give the argument for $p\ge 3$ to simplify the discussion since the argument for $p=2$ is the same up to a change of symbols. 
Recall that the classes $\mu_{n+1}$ and $\lambda_{n+1}$  in $\THH_{*}(\Bn;H\mathbb{F}_p)$ correspond to the comodule primitives $\sigma \btau_{n+1}-\btau_0\sigma \bar{\xi}_{n+1}$ and $\sigma \bar{\xi}_{n+1}$ in $H_{*}(\THH(\Bn;H\mathbb{F}_p))$. 
We therefore have to show that $\sigma \btau_{n+1}-\btau_0\sigma \bar{\xi}_{n+1}$ maps to  $\sigma \bar{\xi}_{n+1}$ under the map $\beta_1$ that is given by applying $H_*(-)$ to 
\[
 \Sigma^{-1} \THH( \Bn;  H\mathbb{F}_p)  \to  \THH( \Bn;  H\mathbb{Z}_{(p)}) \to  \THH( \Bn;  H\mathbb{F}_p). 
   \] 
   As above, one sees that
   \[ H_*(\THH(\Bn;  H\mathbb{Z}_{(p)})) \cong H_*(H\mathbb{Z}_{(p)}) \otimes E(\sigma \bar{\xi}_1, \dots, \sigma\bar{\xi}_{n+1}) \otimes P(\sigma \bar{\tau}_{n+1}). \] 
   The map
   \[ H_*(\THH(\Bn;  H\mathbb{Z}_{(p)}))  \to H_*(\THH(\Bn;  H\mathbb{F}_p))\] 
   is induced by the inclusion $H_*(H\mathbb{Z}_{(p)}) \to  H_*(H\mathbb{F}_p)$. 
   Since the elements $\sigma \bar{\xi}_{n+1}$ and $\sigma \bar{\tau}_{n+1}$ are in the image of this map, they map to zero under $\beta_1$.  Since $\bar{\tau}_0$ is not in the image, it maps to $1$  under $\beta_1$ (up to a unit). Since $\beta_1$ is a derivation, we get $\beta_1(\sigma \btau_{n+1}-\btau_0\sigma \bar{\xi}_{n+1}) \doteq  \sigma \bar{\xi}_{n+1}$.\footnotemark \footnotetext{Note that the Bockstein operator $\beta_1$ is defined for any $H\mathbb{Z}$-algebra $R$ and it is a derivation at this level of generality by \cites{Browder,Shipley}.} Finally, we observe that the $d_1$-differential satisfies the Leibniz rule because the Hurewicz map is a ring map and the Bockstein operator $\beta_1$ is a derivation. 
\end{proof} 

	To compute the differentials $d_r$ for $r > 1$ we use  \cite[Proposition 6.8]{May70}. 
	
\begin{lem}[Proposition 6.8 \cite{May70}]\label{lem: diff in Bockstein ss using may thm}
If $d_{r-1}(x)\ne 0$ in the $v_0$-Bockstein spectral sequence \eqref{v_0Bock} and $|x|=2q$, then 
	\[
	d_{r}(x^p) \dot{=} v_0x^{p-1}d_{r-1}(x)
	\]
	if  $r>2$. If $r=2$ and $p=2$, then 
	\[
	d_{r}(x^p) \dot{=} v_0x^{p-1}d_{r-1}(x)+Q^{|x|}(d_1(x)).
	\]
	If $r=2$ and $p>2$, then 
	\[
	d_{r}(x^p) \dot{=} v_0x^{p-1}d_{r-1}(x)+ \mathsf{E} 
	\] 
	where 
	\begin{equation}\label{eq: obstruction}
	\mathsf{E}= \sum_{j=1}^{(p-1)/2}j[d_1(x)x^{j-1},d_1(x)x^{p-j-1}]_1
	\end{equation}
and $[-,-]_1$ denotes the Browder bracket. 
\end{lem}

\begin{rem}
The result above also appears in \cite{Bru77} in the context of the Adams spectral sequence for an $H_{\infty}$ ring spectrum (cf. \cite[Chapter VI Theorem 1.1,1.2]{BMMS86}) . 
\end{rem}
	
%	The higher Bockstein $\beta_r(x^p)$, which is defined whenever $\beta_{r-1}(x)$ is defined for $ x \in \THH_{2q}(\Bn; H\mathbb{F}_p)$, satisfies: 
%\[ \beta_r(x^p)= 
%\begin{cases}
%	x^{p-1}\beta_1x+Q^{2q}(\beta_1 x) & \text{ if } r=2 , p=2\\
%	x^{p-1} \beta_{r-1} x & \text{ if } r>2 \\
%	x^{p-1} \beta_{r-1} x+ \mathsf{E} & \text{ if } r=2 \text{ and } p>2 
%\end{cases}
%\]
%where
%\begin{align}\label{eq: obstruction} 
%	\mathsf{E}=\sum_{j=1}^{(p-1)/2}j [ \beta_1(x)x^{j-1},\beta_1(x)x^{p-j-1} ]_1.
%\end{align}

%By \cite[Chapter III Theorem 3.2 (1)]{BMMS86}, the Browder bracket $[-,-]_1$ vanishes when we additionally assume that $\Bn$ is an $E_4$ ring spectrum so that $\THH(\Bn; H\mathbb{Z}_p)$ is an $E_3$ ring spectrum. 

We note that in order to apply \cite[Proposition 6.8]{May70}, we need the $\cup_1$-product on $\THH(\Bn;  H\mathbb{Z}_{(p)})$ to satisfy the Hirsch formula, which states that $-\cup_1 c$ is a derivation. We observe that the $\cup_1$-product is a chain homotopy from $x\cdot y$ to $y\cof x$, which corresponds to a braiding in a braided monoidal category. From this perspective, the Hirsch formula corresponds to the first Hexagon axiom in the definition of a braided monoidal category  \cite[\S 1 B1]{JS86}. %[André Joyal and Ross Street, Braided monoidal categories, Macquarie Math Reports 860081 (1986)]. 
It is well documented that there is an $E_2$-operad in small categories with the property that algebras over this operad are braided monoidal categories \cite{Dun97}. The $n$-th category in this operad is the translation groupoid $\text{Br}_n\int \Sigma_n$ of the action of the pure Artin braid group $\text{Br}_n$ on $\Sigma_n$ via the canonical inclusion $\text{Br}_n\to \Sigma_n$. 
We consider the corresponding operad $\mathcal{B}_2$ in $H\mathbb{Z}$-modules by applying the nerve of the category $\text{Br}_n\int \Sigma_n$ and then applying the functor $H\mathbb{Z}_{(p)}\wedge -$. In other words, the $n$-th chain complex in the operad in chain complexes is $\mathcal{B}_2(n)=H\mathbb{Z}_{(p)}\wedge N(\text{Br}_n\int\Sigma_n)_+$. 
The fact that $\THH(B\langle n\rangle ;H\mathbb{Z}_{(p)})$ satisfies the Hirsch formula now follows from two facts: 
\begin{enumerate}
\item algebras over the operad $\mathcal{B}_2$ in chain complexes satisfy the Hirsch formula (cf. \cite[Theorem 1.6]{Dun97}), and 
\item using \cite[Construction 9.6]{May72}, we replace the $E_2$-$H\mathbb{Z}_{(p)}$-algebra $\THH(B\langle n\rangle ;\mathbb{Z}_{(p)})$ with an $\mathcal{B}_2$ algebra without changing the underlying spectrum. 
\end{enumerate}
We therefore tacitly replace our $E_2$ ring spectrum $\THH(B\langle n\rangle ;H\mathbb{Z}_{(p)})$ in $H\mathbb{Z}_{(p)}$-modules with an algebra over the operad $\mathcal{B}_2$ throughout the remainder of the section. The authors thank T. Lawson for suggesting this argument. 

We can consequently prove the following differential pattern. 
\begin{cor}
In the spectral sequence \eqref{v_0Bock}, there are differentials
	\begin{align}\label{primary differentials mod p bock}
	d_{r+1}(\mu_{n+1}^{p^r}) \dot{=} v_0^{r+1}\mu_{n+1}^{p^r-1}\lambda_{n+1}
	\end{align}
	when $p=2$ under the assumption that $B\langle  n\rangle$ is an $E_3$ form. 
	Consequently, there are differentials
	\[
	d_{\nu_p(k)+1}(\mu_{n+1}^k) \dot{=}v_0^{\nu_p(k)+1}\mu_{n+1}^{k-1}\lambda_{n+1}
	\]
	where $\nu_p(k)$ denotes the $p$-adic valuation of $k$.
	The same formulas hold for $p\ge 3$ when the error term \eqref{eq: obstruction} 
	vanishes, for example when $\Bn$ is an $E_{4}$ form of $\tBP{n}$.
	\end{cor}
\begin{proof}
There is a differential
\[
d_1(\mu_{n+1})\dot{=}v_0\lambda_{n+1}
\]
by Lemma \ref{lem: d1 mod p Bockstein} for any prime $p$. 
We will argue that this differential implies the differentials \eqref{primary differentials mod p bock}
for $r\ge 1$ by applying Lemma \ref{lem: diff in Bockstein ss using may thm} and observing that the obstructions vanish. 

When $r=1$ and $p>2$ the the formula \eqref{primary differentials mod p bock} holds whenever the error term \eqref{eq: obstruction} vanishes by Lemma \ref{lem: diff in Bockstein ss using may thm}. The Browder bracket $[-,-]_1$ vanishes by \cite[Proposition 6.3 (iii)]{May70} when $B\langle n \rangle$ is an $E_4$ form of $BP\langle n\rangle$ since in that case 
$\THH(B\langle n\rangle;H\mathbb{Z}_{(p)})$ is an $E_3$ ring spectrum. This completes the base step in the induction for $p>2$.

In the case $p=2$ and $r=1$, Lemma \ref{lem: d1 mod p Bockstein} implies that the error term for $d_2(\mu_{n+1}^2)$ is $Q^{2^{n+2}}\lambda_{n+1}$. At $p=2$, we have that
\begin{align}\label{eq: error term}
Q^{2^{n+2}}\lambda_{n+1} = Q^{2^{n+2}}(\sigma\bxi_{n+1}^2) = \sigma(Q^{2^{n+2}}(\bxi_{n+1}^2)) = \sigma((Q^{2^{n+1}}\bxi_{n+1})^2) = \sigma(\bxi_{n+2}^2)=0
\end{align}
as we now explain.
First, the operation $Q^{2^{n+2}}$ is defined on $\lambda_{n+1}$ because $2^{n+2}=|\lambda_{n+1}|+1$ and $B\langle n\rangle$ is an $E_3$ form of $BP\langle n\rangle$ by assumption. 
The first equality in \eqref{eq: error term} holds by definition of $\lambda_3$, the second equality holds because $\sigma$ commutes with Dyer--Lashof operations by \cite[Proposition 5.9]{AngeltveitRognes}, the third equality holds by \cite[Chapter III Theorem 2.2]{BMMS86}, and the last equality holds because $\sigma$ is a derivation in mod $p$ homology, by \cite[Proposition 5.10]{AngeltveitRognes}. This completes the base step in the induction at $p=2$.

Now let $\alpha=\nu_p(k)$ and let $p$ be any prime. We have that $k=p^\alpha j$ where $p$ does not divide $j$. So by the Leibniz rule
\begin{align*}
		d_{\alpha+1}(\mu_{n+1}^{k}) = & d_{\alpha+1}((\mu_{n+1}^{p^\alpha})^j) \\
		= & j\mu_{n+1}^{p^{\alpha}(j-1)}d_{\alpha+1}(\mu_{n+1}^{p^{\alpha}}) \\
		= & jv_0^{\alpha+1}\mu_{n+1}^{p^\alpha (j-1)}\mu_{n+1}^{p^{\alpha}-1}\lambda_{n+1} \\
		= & v_0^{\alpha+1}\mu_{n+1}^{k-1}\lambda_{n+1}
\end{align*}
	since $j$ is not divisible by $p$ and therefore is a unit in $\mathbb{F}_p$.  
\end{proof}
We now argue that the classes $\lambda_i$ for $1\le i\le n$ are non $p$-torsion in $\THH(\Bn;H\mathbb{Z}_{(p)})$. 
Recall from Proposition \ref{prop: rational thh bpn} that there is an isomorphism
\[\THH_*(\Bn;H\Q)\cong E_{\Q}(\sigma v_1,\dots, \sigma v_{n}).\] 
We claim that the map 
\[\THH_*(\Bn;H\Z_{(p)})\to \THH_*(\Bn;H\Q)\] 
sends $\lambda_i$ to $p^{-1}\sigma v_i$ $1\le i \le n$. To see this, we note that there is a map of $E_2$ ring spectra $\BP\to  \Bn$ by Proposition \ref{properties of forms} and this produces a commutative diagram 
\[ 
	\xymatrix{
	\THH(\BP) \ar[d] \ar[r] & \THH(\BP;H\mathbb{Q}) \ar[d] \\
	\THH(\Bn ;H\mathbb{Z}_{(p)}) \ar[r] & \THH(\Bn ;H\mathbb{Q}) 
	}
\]
of $E_1$ ring spectra by \cite{BFV}. By Proposition \ref{prop: rational thh bpn}, we know $\sigma v_i$ maps to $\sigma v_i$ for $1\le i\le n$ under the left vertical map.
By \cite[Theorem 1.1]{Rog19}, we know that 
\begin{align*}
	\sigma v_i\equiv&p\tilde{\lambda}_i \mod (v_i : i\ge 1) 
\end{align*} 
up to a unit for some classes $\tilde{\lambda}_i=\sigma t_i$. Note that the choice of generators $v_i$ in  \cite[Theorem 1.1]{Rog19} differ from ours, but they are the same up to a unit and modulo decomposables. Therefore there isn't a difference up to a unit modulo $(v_i : i\ge 1)$  after applying the derivation $\sigma$. 
There is an isomorphism 
\[ \THH_*(\BP)\cong E_{\BP_*}(\tilde{\lambda}_k : k\ge 1).\]
and we know that $\tilde{\lambda}_i$ maps to $\lambda_i$ under the map 
\[ \THH_*(\BP)\to \THH_*(\Bn ;\mathbb{Z}_{(p)})\]
for $1\le i\le n$ by Zahler \cite{Zahler} and this does not depend on our choice of $E_3$ form of $\tBP{n}$. 
Therefore, the elements $\lambda_1,\dots ,\lambda_{n}$ are non $p$-torsion and there are no further differentials in the $v_0$-Bockstein spectral sequence \eqref{v_0Bock}.  We define 
\begin{align}\label{HZgens}
		\lambda_{s}:=& 
		\begin{cases} 
		     \lambda_s & \text{ if } 1\le s\le n+1 \\
		\lambda_{s-1}\mu_{n+1}^{p^{s-(n+2)}(p-1)} & \text{ if }s>n+1.
		\end{cases}
\end{align} 
Note that $\THH_*(\Bn;\mathbb{Z}_{(p)})$ is finite type so we can compute $\THH_*(\Bn;\mathbb{Z}_{(p)})$ from $\THH_*(\Bn;\mathbb{Q})$ and $\THH_*(\Bn;\mathbb{Q}_p)$ using the arithmetic fracture square
\[
	\xymatrix{
	\THH_*(\Bn;\mathbb{Z}_{(p)}) \ar[r] \ar[d] & \prod_{p} \THH_*(\Bn;\mathbb{Z}_{p})\ar[d] \\ 
	\THH_*(\Bn;\mathbb{Q}) \ar[r] & \prod_{p}\THH_*(\Bn;\mathbb{Q}_p). 
	}
\]
This proves the following theorem.
\begin{thm}\label{BHZ}
Let $\Bn$ be an arbitrary $E_3$ form of $\tBP{n}$ and at $p>2$ assume the error term \eqref{eq: obstruction} 
vanishes. Then there is an isomorphism of graded $\Z_{(p)}$-modules
\[ \THH_*(\Bn ;H\Z_{(p)})\cong E_{\Z_{(p)}}(\lambda_1,\ldots,\lambda_{n}) \otimes  \left (\Z_{(p)}\oplus T_0^n\right )\]
where $T_0^n$ is a torsion $\Z_{(p)}$-module defined by 
\begin{align}\label{eq: t0n} 
T_0^n=\bigoplus_{s\ge 1}\mathbb{Z}/p^{s}\otimes P_{\mathbb{Z}_{(p)}}(\mu_{n+1}^{p^s}) \otimes \Z_{(p)}\{\lambda_{n+s}\mu_{n+1}^{jp^{s-1}} : 0\le j\le p-2\} 
\end{align}
\end{thm}

%% file: THHBP2coeffHQv3.tex
% root file is THHBP2BP1.tex
\subsection{Rational topological Hochschild homology}
We use the $H\Q$-based B\"okstedt spectral sequence to compute
\[ \pi_*(L_{0}\THH(\Bn))=H\mathbb{Q}_*\THH(\Bn) = \pi_*\THH(\Bn) \otimes \mathbb{Q}\]
for $0\le n\le \infty$ where $B\langle \infty \rangle =\BP$
and $L_0=L_{H\mathbb{Q}}$ is the Bousfield localization at $H\mathbb{Q}$. Since $\BP$ and $\Bn $ are $E_3$ ring spectra, the $H\mathbb{Q}$-based B\"okstedt spectral sequences are strongly convergent multiplicative spectral sequence with signature
\begin{align*} 
E^2_{**}=\HH_{*,*}^{\mathbb{Q}}(H\mathbb{Q}_*\Bn)\implies &H\mathbb{Q}_*\THH(\Bn)
\end{align*} 
for $0\le n\le \infty$. Recall that the rational homology of $\Bn$ is 
\begin{align*}
H\mathbb{Q}_*\Bn\cong &P_\mathbb{Q} (v_1,\dots ,v_n)
\end{align*}
with $|v_i|=2p^i-2$ for $1 \leq i \leq n \leq \infty$. 
Thus, the $E^2$-term of the B\"okstedt spectral sequence is 
\begin{align*}
E^2_{*,*}=& P_\Q(v_1, \dots ,v_n)\otimes_{\mathbb{Q}} E_\Q( \sigma v_1,\dots , \sigma v_n)
\end{align*}
where the bidegree of $\sigma v_i$ is $(1,2(p^i-1))$ for $1 \leq i \leq n \leq \infty$.
Since the $E^2$-page is generated as a $\mathbb{Q}$-algebra by classes in B\"okstedt filtration degree $0$ and $1$, the first quadrant spectral sequence collapses at the $E^2$-page and $E^2_{*,*} = E^\infty_{*,*}$. 
There are no  multiplicative extensions, because the $E^{\infty}$-pages are free graded-commutative $\mathbb{Q}$-algebras. 
Therefore, we produce isomorphisms of graded $\Q$-algebras  
\begin{align*}
\THH_*(\Bn)\otimes \Q\cong & P_\Q(v_1,\dots,v_n)\otimes_\Q E_\Q(\sigma v_1,\dots,\sigma v_n) 
\end{align*}
with $|\sigma v_i|=2p^i-1$ for $1 \leq i \leq n \leq \infty$. It follows that there is an equivalence 
\[L_0 \THH(\Bn)\simeq \bigvee_{x\in B_n} \Sigma^{|x|}L_0\Bn,\]
where $B_n$ is a graded basis for $E_\Q(\sigma v_1,\dots,\sigma v_n)$ as a graded $\mathbb{Q}$-vector space, since $L_0$ is a smashing localization. We may also let $n=\infty$ and in this case $B\langle \infty \rangle=\BP$ and $B_{\infty}$ is a graded basis for $E_\Q(\sigma v_1,\sigma v_2,\dots)$ as a graded $\mathbb{Q}$-vector space. 

By Proposition \ref{properties of forms}, the linearization map $\tBP{n}\to H\mathbb{Z}_{(p)}$ is an $E_3$ ring spectrum map. Since the localization map $H\mathbb{Z}_{(p)}\to H\mathbb{Q}$ is an $E_{\infty}$ ring spectrum map, we may infer that the B\"okstedt spectral sequence 
\begin{align*} 
E^2_{**}=\HH_{*,*}^{\mathbb{Q}}(H\mathbb{Q}_*\Bn;\mathbb{Q})\implies &H\mathbb{Q}_*\THH(\Bn; H\mathbb{Q})
\end{align*} 
is a spectral sequence of $\mathbb{Q}$-algebras by adapting the proof of   \cite[Proposition 4.3]{AngeltveitRognes} using \cite[\S 3.3]{BFV}. 
This spectral sequence collapses without extensions by the same argument as before. All of these computations are functorial with respect to the map of $E_2$-ring spectra $\BP\rightarrow \Bn$ from Proposition \ref{properties of forms}. This proves the following result. 
\begin{prop}\label{prop: rational thh bpn}
There is an isomorphism of graded $\mathbb{Q}$-algebras
\begin{equation}\label{Qcoeff} \THH_*(\Bn;H\mathbb{Q})\cong E_{\mathbb{Q}}(\sigma v_1,\dots , \sigma v_n)
\end{equation}
for all $0 \le n \le \infty$. The map 
\[ \THH_*(\BP; H\mathbb{Q})\to \THH_*(\Bn ; H\mathbb{Q}) \]
sends $\sigma v_i$ to $\sigma v_i$ for $0\le i\le n$. 
\end{prop}

%% file: THHBP2coeffk1v3.tex
% root file is THHBP2BP1.tex
\section{Topological Hochschild homology mod $(p,v_2)$}\label{sec: THH BP2 mod p, v2}

In this section, we compute topological Hochschild homology of $\Bt$ with coefficients in $k(1)$. First we compute topological Hochschild homology with coefficients in $K(1)$.

\input{THHBP2K1v3}

\subsection{The $v_1$-Bockstein spectral sequence}\label{subsec: v1BSS}
We compute $\THH_*(\Bt; k(1))$ using the spectral sequence \eqref{Bockstein vi} for $n=2$ and $i=1$.
For $s \geq 4$, we recursively define 
\[ \lambda_s :=  \lambda_{s-2} \mu_3^{p^{s-4}(p-1)}.\]  
For $s \geq 1$, we  define 
\[ r(s,1) := \begin{cases} p^{s+1} + p^{s-1} + \dots + p^2 & s \equiv 1 \; \text{mod} \; 2 \\
p^{s+1} + p^{s-1} + \dots + p^3  & s \equiv  0 \; \text{mod} \;  2. 
\end{cases}  \] 

\begin{thm} \label{mod p v_2}
Let $\Bt$ be an $E_3$ form of $BP \langle 2 \rangle$ and let $p \geq 3$. There is an isomorphism of $P(v_1)$-modules 
\[\THH_*(\Bt; k(1)) \cong E(\lambda_1) \otimes (P(v_1) \oplus T^2_1), \] 
where 
\begin{align}\label{eq: t12} T^2_1 = \bigoplus_{s \geq 1} P_{r(s,1)}(v_1) \otimes E(\lambda_{s+2}) \otimes P(\mu_3^{p^s}) \otimes \mathbb{F}_p\{\lambda_{s+1} \mu_3^{jp^{s-1}}: 0 \leq j \leq p-2\}.
\end{align}	
\end{thm} 
\begin{proof} 
We prove by induction on $s \geq 1$ that 
\[ E_{r(s,1)}^{*,*} = E(\lambda_1) \otimes  \bigl(P(v_1) \otimes E(\lambda_{s+1}, \lambda_{s+2}) \otimes P(\mu_3^{p^{s-1}})  \oplus M_s\bigr) \] 
with
\[ M_s = \bigoplus_{t =1}^{s-1} P_{r(t,1)}(v_1) \otimes E(\lambda_{t+2}) \otimes P(\mu_3^{p^t}) \otimes \mathbb{F}_p\{\lambda_{t+1}\mu_3^{jp^{t-1}} : 0 \leq j \leq p-2\},  \] 
that we have a differential 
$d_{r(s,1)}(\mu_3^{p^{s-1}}) \dot{=} v_1^{r(s,1)} \lambda_{s+1}$,   and that the classes $\lambda_{s+1}$ and $\lambda_{s+2}$ are infinite cycles. This implies the statement. 

By Theorem \ref{thm:K(1)coeff}, the elements $v_1^s$ are permanent cycles for every $s$, so the classes $\lambda_2$ and $\lambda_3$ cannot support differentials and thus are infinite cycles.  Note that we use $p \geq 3$ here; for $p = 2$ we  would  have a possible differential $d_2(\lambda_3) \dot{=} v_1^2\lambda_1\lambda_2$. Since the classes $v_1^n\lambda_1$ survive by Theorem \ref{thm:K(1)coeff}, the only possible differential on $\mu_3$ is 
\[ d_{p^2}(\mu_3) \dot{=} v_1^{p^2}\lambda_2\] 
for bidegree reasons.
This differential must exist because otherwise the spectral 
sequence  would collapse at the $E_2$-page by multiplicativity  which would contradict Theorem \ref{thm:K(1)coeff}.   This proves the base step $s = 1$ of the induction.  Now, assume that the statement holds for some $s \geq 1$. 
We then get 
\[ E_{r(s,1)+1} ^{*,*} = E(\lambda_1) \otimes (P(v_1) \otimes E(\lambda_{s+2},\lambda_{s+1}\mu_3^{p^{s-1}(p-1)} ) \otimes P(\mu_3^{p^{s}}) \oplus M_{s+1},\] 
and it suffices to show that
$\lambda_{s+3} = \lambda_{s+1}\mu_3^{p^{s-1}(p-1)}$ is an infinite cycle and that 
\[ d_{r(s+1,1)}(\mu_3^{p^s}) \dot{=} v_1^{r(s+1,1)}\lambda_{s+2}.\] 
Note that the class $\lambda_{s+2}$ is an infinite cycle by the induction hypothesis. The class $\lambda_{s+3}$ is an infinite cycle for bidegree reasons and because the classes $v_1^s$ are permanent cycles. Note that we use $p \geq 3$ here; for $p = 2$ and $s$ even we would have a possible differential 
$d_{r(s,1)+p}(\lambda_{s+3}) \dot{=}v_1^{r(s,1)+p} \lambda_1 \lambda_{s+2}$. 
The class $\mu_3^{p^s}$ must support a differential because otherwise the spectral sequence would collapse at this stage which would contradict Theorem \ref{thm:K(1)coeff}.
Since the classes $v_1^n \lambda_1$ are permanent cycles,
 we get 
 \[ d_{r(s+1,1)}(\mu_3^{p^s}) \dot{=} v_1^{r(s+1,1)}\lambda_{s+2}\] 
 for bidegree reasons. 
 Here note that $v^{r(s,1)}_1\lambda_{s+3}$ has the right topological degree, but the filtration degree is too low for it to be the target of a differential on $\mu_3^{p^s}$ at the $E_{\ell}$-page for $\ell>r(s,1)$. This completes the induction step. 
\end{proof}

%% file: THHBP2K1v3.tex
\subsection{$K(1)$-local topological Hochschild homology }\label{coefficientsK1}
In this section we assume that $p \geq 3$ and we write $\Bt$ for an $E_{3}$ form of $\B$. Write $k(1)=\Bt/(p,v_2)$ for the $E_1$ $\Bt$-algebra constructed as in Proposition \ref{properties of forms} and let $K(1)=k(1)[v_1^{-1}]$. In order to determine the topological Hochschild homology of $\Bt$ with coefficients in $k(1)$, we first determine 
\[\THH(\Bt; K(1))=THH(\Bt; K(1)).\]   
To compute the multiplicative B\"okstedt spectral sequence
\[ E^2_{*,*} = \HH_{*,*}^{K(1)_*}(K(1)_*\Bt) \Longrightarrow K(1)_*\THH(\Bt)\] 
we first need to compute $K(1)_*\Bt$. 
To compute $K(1)_*\Bt$ we first relate it to $\BP_*\BP$. Recall that we have 
\[ \BP_*\BP = BP_*[t_1, t_2, \dots] \] 
with $|t_i| = 2p^i-2$. By \cite[Theorem A2.2.6]{greenbook} the right unit $\eta_R$ is determined by
\begin{equation} \label{eq: right unit in BP coops}
 \sum^F_{i,j \geq 0} t_i \eta_R(v_j)^{p^i} = \sum^F_{i,j \geq 0} v_i t_j^{p^i}, 
 \end{equation}
 where $t_0 = 1$ and $v_0 = p$. 
\begin{lem} \label{tens}
The composite map 
 \[ K(1)_* \otimes_{\BP_*} \BP_*\BP \otimes_{\BP_*} \Bt_* \to \pi_*(K(1) \wedge_{\BP} (\BP \wedge \BP) \wedge_{\BP}  \Bt) \cong K(1)_*\Bt\] 
is an isomorphism. 
\end{lem}
\begin{proof}
Consider the commutative diagram 
\begin{equation} \label{tensorsp}
\begin{tikzcd}
\pi_*(K(1) \wedge \Bt) \ar{r} & \pi_*(K(1) \wedge \Bt[v_1^{-1}]) \\
\pi_*(K(1) \wedge_{\BP} (\BP \wedge \BP) \wedge_{\BP}  \Bt) \ar{r} \ar{u}{\cong} & \pi_*(K(1) \wedge_{\BP} (\BP \wedge \BP) \wedge_{\BP} \Bt[v_1^{-1}] ) \ar{u}{\cong}\\
K(1)_* \otimes_{\BP_*} \BP_*\BP \otimes_{\BP_*} \Bt_* \ar{r} \ar{u}  & \pi_*(K(1)) \otimes_{\BP_*} \BP_*\BP \otimes_{\BP_*} \Bt_*[v_1^{-1}]  \ar{u}. 
\end{tikzcd}
\end{equation}
Since $\Bt[v_1^{-1}]$ is Landweber exact, the right-hand  vertical  map is an isomorphism.  
In (\ref{eq: right unit in BP coops}) the $F$-summands in degree $\leq 2p-2$
are $\eta_R(v_0)$, $t_1\eta_R(v_0)^p$, $\eta_R(v_1)$, $ v_0$,  $v_1$ and $v_0t_1$. 
Thus, we have $\eta_R(v_1) = v_1$ in $K(1)_* \otimes_{\BP_*} \BP_*\BP = K(1)_*[t_i \mid i\ge 1]$, because we have  $p = 0$ in this ring. We  get that in $K(1)_* \otimes_{\BP_*} \BP_*\BP \otimes_{\BP_*} \Bt_* $
\[ v_1 \otimes 1 \otimes 1 = 1 \otimes v_1 \otimes 1 = 1 \otimes \eta_R(v_1) \otimes 1 = 1 \otimes 1 \otimes v_1\] 
holds. This implies that the upper and lower horizontal map in the diagram are isomorphisms. It follows that the left vertical map is an isomorphism too. 
\end{proof}

\begin{not*}
Let $f_i(v_1, v_2)  \in \Bt_* = \mathbb{Z}_{(p)}[v_1, v_2]$ be  the image of $v_i$ under $\BP_* \to \Bt_*$. Define 
\[ v_i' \coloneqq v_i - f_i(v_1, v_2)  \in \BP_*. \]  Then, $v_i'$ is in the kernel of $\BP_* \to \Bt_*$ and $\BP_* = \mathbb{Z}_{(p)}[v_1, v_2, v_3', \dots]$.
\end{not*}

 By Lemma \ref{tens} we get
\begin{align*}
 K(1)_*\Bt & =   (K(1)_* \otimes_{BP_*} \BP_*[t_1, \dots]) \otimes_{\mathbb{Z}_{(p)}[v_1, v_2, v_3^{\prime} 
  \dots]} \mathbb{Z}_{(p)}[v_1, v_2] \\
&  =     K(1)_*[t_i\mid i\ge 1]/{(\eta_R(v_3'), \dots)}.
  \end{align*}
  
 \begin{lem}  \label{lem: ci}
 For $i \geq 0$ the element $\eta_R(v_{i+1}) \in K(1)_*[t_i\mid i\ge 1]$  actually lies in $K(1)_*[t_1, \dots, t_{i}]$. In fact, we have 
 \[
 	\eta_R(v_{i+1}) = v_{i+1} + v_1t_i^p - v_1^{p^{i}}t_i + g_{i}
 \]
 where $g_{i}\in K(1)_*[t_1, \ldots, t_{i-1}]$. 
  \end{lem}
  \begin{proof}
  We will prove the claim in 
  $k(1)_*[t_i\mid i\ge 1]$, from this the result will follow. The reason we do this is because we will want to make degree arguments, and hence will want to avoid negative gradings. 
  
  In $BP_*BP/(p)$, we have $\eta_R(v_1)  =v_1$. It also follows from \eqref{eq: right unit in BP coops} that for $i\geq 0$
  \[
  	\eta_R(v_{i+1}) \equiv v_{i+1} + v_1t_i^p-v_1^{p^i}t_i \mod (t_1, t_2, \ldots, t_{i-1})
  \] 
  in $BP_*BP/(p)$. Thus, this congruence also holds in $k(1)_*[t_i \mid i\ge 1]$. Since $\eta_R(v_{i+1})$ lifts to $BP_*BP/(p)$ we may make our degree arguments in $k(1)_*[t_i \mid i\ge 1]$. In the ring $k(1)_*[t_i \mid i\ge 1]$, we therefore have that 
  \[
  	\eta_R(v_{i+1}) = v_{i+1} + v_1t_i^p-v_1^{p^i}t_i + g_i
  \]
  where $g_i$ is a polynomial in the ideal generated by $t_1, t_2, \ldots, t_{i-1}$. Thus far we have not excluded the possibility that a monomial divisible by  $t_j$ with $j\geq i$ occurs as a summand of $g_i$. 
  
  For $j>i+1$, we can exclude this possibility for degree reasons. Indeed, $\eta_R(v_{i+1})$ is homogenous of degree $2(p^{i+1}-1)$, and when $j>i+1$ the element $t_j$ has degree greater than $2(p^{i+1}-1)$. Consider the case when $j=i+1$. To exclude this case, suppose there exists a monomial $m$ in $k(1)_*[t_1, \ldots,t_i]$ which is a summand of $g_i$ and is divisible by $t_{i+1}$. Then as the degrees of $t_{i+1}$ and $\eta_R(v_{i+1})$ are the same, it follows that $m = at_{i+1}$ for some $a\in \F_p$. If $a\neq 0$, then this contradicts the assumption that $g_i$ is in the ideal $(t_1, \ldots, t_{i-1})$. This shows that $g_i\in k(1)_*[t_1, \ldots, t_{i}]$. 
  
  We now exclude the possibility that a monomial divisible by $t_i$ occurs as a summand of $g_i$. Note that the summands $t_k\eta_R(v_j)^{p^k}$ and $v_kt_j^{p^k}$ in \eqref{eq: right unit in BP coops} both have degree $2(p^{k+j}-1)$. Cross terms in \eqref{eq: right unit in BP coops} from those summands with degree less than or equal than $2(p^{i+1}-1)$ could potentially produce a $t_i$ divisible monomial as a summand of $g_i$. On the right-hand side of \eqref{eq: right unit in BP coops}, the possible summands are those of the form $v_jt_i^{p^j}$. As this must have degree at most $2(p^{i+1}-1)$, we must have $j=0,1$. These correspond, respectively, to $v_0t_i=pt_i$ and $v_1t_i^p$. But $p=0$ in $k(1)_*$, so the only one to consider is $v_1t_i^p$, This has degree exactly $2(p^{i+1}-1)$, and so a monomial divisible by this element does not occur in $g_i$. In fact, it has already been accounted for. 

  For the left-hand side, we similarly find that the only summand which could potentially produce a $t_i$ divisible monomial as a summand of $\eta_R(v_{i+1})$ is
  \[
  	t_i\eta_R(v_1)^{p^i} = t_iv_1^{p^i}.
  \]
  As this has exactly degree $2(p^{i+1}-1)$, it does not occur in $g_i$ because it cannot be written as an element in the ideal $(t_1,\dots t_{i-1})$ for degree reasons.
  In fact, this element has already been accounted for. Thus there are no $t_i$ divisible monomials appearing as summands of $g_i$. Consequently, we have shown that $g_i\in k(1)_*[t_1, \ldots, t_{i-1}]$ as desired. 
  \end{proof}
  Recall from the proof of Proposition \ref{properties of forms} that we have 
  \[
  	v_i' = v_i-f_i(v_1,v_2)
  \]
  for some $f_i\in \mathbb{Z}_{(p)}[x,y]$. In light of the previous lemma, we conclude that the class 
  \[
  	\eta_R(v_i') = \eta_R(v_i) - f_i(\eta_R(v_1), \eta_R(v_2)) \in K(1)_*[t_i\mid i\ge 1]
  \] 
  also lies in $K(1)_*[t_1, \dots, t_{i-1}]$ for each $i \geq 3$. 
 \begin{lem} 
 The maps of commutative $K(1)_*$-algebras
 \[ K(1)_*[t_1, \dots, t_{i-1}] /({\eta_R(v_3^{\prime}), \dots, \eta_R(v_i^{\prime})}) \to K(1)_*[t_1, \dots, t_{i}] /({\eta_R(v_3^{\prime}), \dots, \eta_R(v_{i+1}^{\prime})})  \] 
induced by precomposing the canonical quotient map with the canonical inclusion map are \'etale for $i \geq 2$. 
 \end{lem}
 \begin{proof}
 	For ease of notation, set 
 	\[
 		A_i:= K(1)_*[t_1, \ldots, t_{i-1}]/(\eta_R(v_3'), \ldots, \eta_R(v_i'))
 	\]
 	for $i\geq 2$. Note that Lemma \ref{lem: ci} allows us to make this definition. Note also that the $A_{i+1} = A_i[t_i]/(\eta_R(v_{i+1}'))$. We wish to show that the map
 	\[
 		A_i\to A_{i+1}
 	\] 
 	is an \'etale morphism. To do this, it is enough to show that the partial derivative of $\eta_R(v_{i+1}')$ with respect to $t_i$ is a unit in $A_i$. Write $\partial_i$ for the partial derivative with respect to $t_i$. 
 	Since
 	\[
 		v_{i+1}' = v_{i+1}-f_{i+1}(v_1,v_2)
 	\]
 	for some $f_{i+1}\in \mathbb{Z}_{(p)}[x,y]$, we can infer that  
 	\[
 		\eta_R(v_{i+1}') = \eta_R(v_{i+1})-f_{i+1}(\eta_R(v_1),\eta_R(v_2)). 
 	\]
 	In $K(1)_*[t_1, t_2,\ldots]$, we know 
 	$
 		\eta_R(v_1) = v_1
 	$
 	since $p=0$ in $K(1)_*$, and we have 
 	\[
 		\eta_R(v_2) = v_1t_1^p-v_1^pt_1.
 	\]	
 	Thus, 
 	\[
 		\partial_i \eta_R(v_{i+1}') = \partial_i \eta_R(v_{i+1}),
 	\]
 	for $i\geq 2$ and it suffices to show that $\partial_i \eta_R(v_{i+1})$ is a unit in $A_i$. 
 	
 	By Lemma \ref{lem: ci}, we have the formula  
 	\[
 		\eta_R(v_{i+1}) = v_{i+1} + v_1 t_i^p- v_1^{p^i} t_i +g_i
 	\]
 	where $g_i \in K(1)_*[t_1, \ldots, t_{i-1}]$. Thus, we conclude that  
 	\[
 		\partial_i \eta_R(v_{i+1})= -v_1^{p^i}\in K(1)_*[t_1, \ldots, t_{i-1}].
 	\]
 	Since $v_1^{p^i}$ is a unit in $K(1)_*$, this shows that $\partial_i\eta_R(v_{i+1})$ is a unit in $A_i$. 
\end{proof}

We continue to use the notation from the proof of the previous lemma. Since each map $A_i\to A_{i+1}$ is \'etale, we may apply  \cite[Theorem 0.1]{WeibelGeller} to conclude that 
\begin{align*}
\HH_{*,*}^{K(1)_*}(K(1)_*\Bt) & =  \colim \HH_{*,*}^{K(1)_*}(A_i) \\  
& =  \colim \HH^{K(1)_*}_{*,*}(A_2) \otimes_{A_2} A_i \\
& =  \HH_{*,*}^{K(1)_*}(A_2) \otimes_{A_2} K(1)_*\Bt \\
& =  E(\sigma t_1) \otimes K(1)_*\Bt. 
\end{align*}
Since this is concentrated in B\"okstedt filtration $0$ and $1$, the B\"okstedt spectral sequence collapses yielding
\[ E(\sigma t_1) \otimes K(1)_*\Bt \cong K(1)_*\THH(\Bt).\] 

In the Hopf algebroid $(\BP_*, \BP_*\BP)$, we have the formula
\[ 
\sum^F_{i \geq 0} \; \Delta(t_i) = \sum^F_{i,j \geq 0} \;  t_i \otimes t_j^{p^i}
\] 
by \cite[Theorem A2.1.27]{greenbook}. Since the $\BP_*\BP$-coaction on $t_i$ agrees with the coproduct, it is determined by the formula 
\[\Delta(t_1) = 1 \otimes t_1 + t_1 \otimes 1.\] 
Note that $(K(1)_*, K(1)_*K(1))$ is a flat Hopf algebroid and $K(1)_*(X)$ 
is a left $K(1)_*K(1)$-comodule for every spectrum $X$. 
By naturality, we observe that $t_1  \in K(1)_*\Bt$  has the $K(1)_*K(1)$-coaction $1 \otimes t_1 +  t_1 \otimes 1$. 
Let
\[ \sigma \colon \thinspace K(1)_*\Bt \to K(1)_{*+1}\THH(\Bt) \] 
be the usual $\sigma$ operator analogous to the one defined in \cite{McClureStaffeldt}.
By \cite[Proposition 5.10]{AngeltveitRognes}, which also applies to our setting because the Hopf element $\eta=0\in K(1)_*$, the operator $\sigma$ is a derivation.
It is also clear that $\sigma$ is compatible with the $K(1)_*K(1)$-comodule action in the sense that 
\[ \psi(\sigma x)=(1\otimes \sigma)(\psi(x))\]
where 
\[ \psi \colon \thinspace K(1)_{*}\THH(\Bt)\to K(1)_*K(1)\otimes K(1)_{*}\THH(\Bt).\]
It follows that $\sigma t_1 \in K(1)_*\THH(\Bt)$ is a comodule primitive.  
Since there is a weak equivalence $\THH(B\langle 2\rangle, K(1)) \simeq K(1) \wedge_{B\langle 2 \rangle } THH(B\langle 2 \rangle)$ by \cite[Remark 6.1.4]{HahnWilson20}, we may infer from the K\"unneth isomorphism that there is an isomorphism of $K(1)_*$-modules
\[ K(1)_*\THH(\Bt; K(1)) \cong  K(1)_*K(1) \otimes E(\sigma t_1),\]
where $\sigma t_1$ is a comodule primitive. 
Since $\THH(\Bt; K(1))$ is a $K(1)$-module spectrum and $K(1)_*$ is a graded field, we have that it splits as a sum of suspensions of $K(1)$ and that its homotopy is isomorphic to the comodule primitives in  $K(1)_*\THH(\Bt; K(1))$. 
 Thus, there is an isomorphism of $K(1)_*$-modules
\[ \THH_*(\Bt; K(1)) = K(1)_* \otimes E(\sigma t_1).\] 
Since $\sigma t_1$ lifts to a class in $\tilde{\lambda}_1\in\THH_*(\Bt; k(1))$ which projects onto $\lambda_1$ via the map 
$\THH_*(\Bt; k(1))\to \THH_*(\Bt; H\mathbb{F}_p)$
induced by the linearization map $k(1)\to H\mathbb{F}_p$ by \cite{Zahler}, we simply rename this class $\lambda_1$. 

In summary, we have proven the following theorem. 
\begin{thm}\label{thm:K(1)coeff} For $\Bt$ an $E_3$ form of $\tBP{2}$ and $p\geq 3$, the following hold:
	\begin{enumerate}
		\item There is a weak equivalence
		\[ K(1)\vee \Sigma^{2p-1}K(1)\simeq \THH(\Bt;K(1)).\]
		\item The $P(v_1)$-module $\THH_*(\Bt;k(1))$, modulo $v_1$-torsion, is freely generated by $1$ and $\lambda_1$.
	\end{enumerate}
\end{thm}

%% file: THHBP2coeffk2v3.tex
% root file is V(0)THHBP2.tex

\section{Topological Hochschild homology mod $(p,v_1)$}\label{sec: THH BP2 mod p, v1}

In this section $\Bt$ is again  an $E_3$ form of $BP \langle 2 \rangle$, e.g. $\tmf_1(3)$ at $p = 2$, $\taf$ at $p = 3$, or $BP\langle n\rangle^{\prime}$ at an arbitrary prime $p$. 
We let $k(2) \coloneqq \Bt/{(p,v_1)}$ be the $E_1$-$\Bt$-algebra contructed in Proposition \ref{properties of forms} and let $K(2)=k(2)[v_2^{-1}]$. The goal of this section is to compute the homotopy groups of $\THH(\Bt; K(2))$. In Subsection \ref{thhbtK2}, we first show that the unit map 
\[ K(2) \longrightarrow \THH_*(\Bt; K(2))\] 
is an equivalence. 
This implies that in the abutment of the $v_2$-Bockstein spectral sequence
 \[ \THH_*(\Bt; H\mathbb{F}_p)[v_2] \Longrightarrow \THH_*(\Bt; k(2)) \] 
 all classes are $v_2$-torsion besides the powers of $v_2$.  This  allows us to compute this spectral sequence in Subsection \ref{BSSv2}.

\subsection{$K(2)$-local topological Hochschild homology} \label{thhbtK2}
Considering a diagram analogous to (\ref{tensorsp}), one sees that we have an isomorphism 
\[\begin{tikzcd}
 K(2)_* \otimes_{\BP_*} \BP_*\BP \otimes_{\BP_*} \Bt_* \ar{r} & \pi_*(K(2) \wedge \Bt).
 \end{tikzcd}\] 
 For this, note that 
 \[\eta_R(v_1) =  v_1 = 0 \in  K(2)_* \otimes_{\BP_*} \BP_*\BP =K(2)_*[t_i\mid i\ge 1]  \] 
 and therefore $\eta_R(v_2) = v_2$. This implies that the equality 
 \[ v_2 \otimes 1 \otimes 1 = 1 \otimes 1 \otimes v_2\]
 holds in the tensor product 
 \[ K(2)_* \otimes_{\BP_*} \BP_*\BP \otimes_{\BP_*} \Bt_*.\] 
From this, we determine that  
\[
	K(2)_*\Bt = K(2)_*[t_i\mid i\ge 1]/{(\eta_R(v_3'), \dots)}.
\]
In particular, this is a graded commutative $K(2)_*$-algebra even at $p=2$ where $K(2)$ is not homotopy commutative (cf. \cite[Lemma 8.9]{AngeltveitRognes}). 

\begin{lem}\label{lem: error term in height 2 right unit}
	In $K(2)_*[t_1 \mid i\ge 1]$ we have that 
	\[
		\eta_R(v_{i+2}) = v_{i+2} + v_2t_i^{p^2} - v_2^{p^i}t_i +g_i
	\]
	where $g_i\in K(2)_*[t_1, \ldots, t_{i-1}]$. 
\end{lem}
\begin{proof}
	We argue similarly to Lemma \ref{lem: ci} and make our arguments in the ring $k(2)_*[t_i \mid  i \ge 1]$. The result will follow from this. We have that 
	\[
	\eta_R(v_{i+2}) \equiv v_{i+2} + v_2t_i^{p^2} - v_2^{p^i}t_i \mod (t_1, t_2, \ldots, t_{i-1}),
	\]
	in $\BP_*\BP/(p,v_1)$ (c.f. \cite[Proof of Theorem 4.3.2]{greenbook}). Consequently, this formula also holds in $k(2)_*[t_i\mid i\geq 1]$. This shows that in $k(2)_*[t_i\mid i\geq 1]$ we have
	\[
		\eta_R(v_{i+2}) = v_{i+2} + v_2t_i^{p^2}-v_2^{p^i}t_i+g_i
	\]
	for some $g_i$ in the ideal $(t_1,t_2, \ldots, t_{i-1})$. 
	Since $\eta_R(v_{i+2})$ lifts to the graded abelian group $BP_*BP/(p,v_1)$, we may also make degree arguments in $k(2)_*[t_i \mid i\ge 1]$. 
	
	Note that for degree reasons, there can be no instance of a $t_j$ with $j>i+2$ dividing a monomial summand of $g_i$. 
	We can also exclude the possibility of $t_{i+2}$ dividing a monomial in $g_i$. 
	Indeed, a monomial in $g_i$ divisible by $t_{i+2}$ would necessarily be just $t_{i+2}$ itself, contradicting that $g_i$ is in the ideal $(t_1, \ldots, t_{i-1})$. 
	This shows that we have 	
	\[
		\eta_R(v_{i+2})\in k(2)_*[t_1, \ldots, t_{i+1}]. 
	\]
	for all $i\geq 0$. 

	We now exclude the possibility that $t_{i+1}$ divides a monomial in $\eta_R(v_{i+2})$. To do this, we note that a $t_{i+1}$ divisible monomial in $g_i$ could arise from cross terms involving the universal $p$-typical formal group law and the formula \eqref{eq: right unit in BP coops}. Note that the only terms to consider on the right hand side are $v_0t_{i+1}$ and $v_1t_{i+1}^p$, which are $0$ since $p=v_1=0\in k(2)_*$. On the left hand side, we only need to consider the terms $t_k\eta_R(v_{j+2})^{p^k}$ of degree less than or equal to $2(p^{i+2}-1)$. This immediately implies that $j\leq i$. For $k=i+1$, the term of smallest degree is $t_{i+1}\eta_R(v_2)^{p^{i+1}}$. The degree of this term is $2(p^{i+3}-1)$, which is too large. Thus we can exclude the possibility that $k=i+1$. Now as $j\leq i$ and since we have shown that $\eta_R(v_{j+2})\in k(2)_*[t_1, \ldots, t_{j+1}]$, we see that none of the relevant terms on the left hand side can contribute a $t_{i+1}$ divisible monomial summand to $\eta_R(v_{i+2})$. Thus we have that $g_i\in K(2)_*[t_1, \ldots, t_i]$.
		
	We are left to consider whether a $t_i$ divisible monomial could occur as a summand of $g_i$ via the cross terms coming from the formal group law $F$ in \eqref{eq: right unit in BP coops}. On the right hand side, we only need to consider the term $v_2t_i^{p^2}$. Here we use the fact that $v_1=0\in k(2)_*$. This term has already been accounted for and is not in $g_i$. On the left hand side, since we have shown that $\eta_R(v_{j+2})\in k(2)_*[t_1, \ldots, t_j]$, the only term we need to consider is $t_iv_2^{p^i}$. Again, we have already considered this term. We can therefore conclude that $g_i\in k(2)_*[t_1, \ldots, t_{i-1}]$.
\end{proof}

 \begin{defn}\label{ai and hi}
 We define commutative $K(2)_*$-algebras 
 \begin{align*}
 C_0 & \coloneqq K(2)_* \\
 C_i & \coloneqq C_{i-1}[t_i]/{\eta_R(v_{i+2}')}, \;\;\; i \geq 1
 \end{align*}
 and write $h_i \colon \thinspace C_{i-1}\rightarrow C_{i}$ for the map of commutative $K(2)_*$-algebras defined as the 
 composite of the canonical inclusion map $C_{i-1}\rightarrow C_{i-1}[t_i]$ with the canonical quotient map $C_{i-1}[t_i]\rightarrow  C_{i-1}[t_i]/{\eta_R(v_{i+2}')}$. 
 \end{defn}
 Thus we have 
 \[
 	C_i = K(2)_*[t_1, \ldots, t_i]/(\eta_R(v_3'), \ldots, \eta_R(v_{i+2}'))
 \]
 for $i\geq 1$ and 
 \[ K(2)_*\Bt = \colim_i C_i.\]
 We proceed in the same fashion as in Section \ref{coefficientsK1} and argue that $h_i\colon \thinspace C_{i-1} \to C_i$ is \'etale by examining the derivative of $\eta_R(v_{i+2}')$ with respect to $t_i$.
 
 \begin{lem}\label{lem: height 2 transition maps are etale}
 	The map of commutative rings $h_i\colon \thinspace C_{i-1}\to C_i$ from Definition \ref{ai and hi} is \'etale. 	
 \end{lem}
 \begin{proof}
 	We have that 
 	\[
 		v_{i+2}' = v_{i+2}-f_{i+2}(v_1,v_2) = v_{i+2}-f_{i+2}(0,v_2). 
 	\]
 	Hence we have that 
 	\[
 		\eta_R(v_{i+2}') = \eta_R(v_{i+2}) - f_{i+2}(0,v_2). 
 	\]
 	Let $\partial_i$ denote the partial derivative with respect to $t_i$. Since $C_i = C_{i-1}[t_i]/(\eta_R(v_{i+2}'))$, to show the morphism $C_{i-1}\to C_i$ is \'etale, it is enough to show that $\partial_i\eta_R(v_{i+2}')$ is a unit. We have 
 	\[
 		\partial_i \eta_R(v_{i+2}') = \partial_i\eta_R(v_{i+2}) - \partial_i f_{i+2}(0,v_2) = \partial_i \eta_R(v_{i+2}). 
 	\]
 	From Lemma \ref{lem: error term in height 2 right unit}, we find that $\partial_i g_i=0$, and hence
 	\[
 		\partial_i \eta_R(v_{i+2}) = \partial_i \left(v_{i+2} + v_2t_i^{p^2} - v_2^{p^i}t_i +g_i\right) = -v_2^{p^i}
 	\]
 	which is a unit. This completes the proof. 
 \end{proof}
 
 Since each map $C_i\to C_{i+1}$ is \'etale, we may apply  \cite[Theorem 0.1]{WeibelGeller} to conclude that the unit map 
\begin{align}\label{algebraic unit map}
K(2)_*\Bt \to \HH^{K(2)_*}_{*,*}(K(2)_*\Bt)
\end{align} 
 is an isomorphism of graded commutative $\mathbb{F}_p$-algebras (even at $p=2$). 
 The unit map $K(2)_*\Bt \to K(2)_*\THH(\Bt)$ is the edge homomorphism in the B\"okstedt spectral sequence
 \[ E^2_{*,*} = \HH^{K(2)_*}_{*,*}(K(2)_*\Bt) \Longrightarrow K(2)_*\THH(\Bt) \]
and the input is concentrated in B\"okstedt filtration zero by \eqref{algebraic unit map}, so the spectral sequence collapses without extensions yielding an isomorphism 
\[ K(2)_*\Bt \cong  K(2)_*\THH(\Bt)\]
of graded commutative $\mathbb{F}_p$-algebras (even at the prime $p=2$). 
 
 By the K\"unneth isomorphism, the map 
\[\begin{tikzcd}
 K(2)_*K(2) \ar{r} & K(2)_*\THH(\Bt, K(2)) 
 \end{tikzcd} \] 
 is an isomorphism as well. Since both $K(2)$ and $\THH(\Bt; K(2))$ are $K(2)$-local, we obtain the following result.
\begin{cor}\label{cor:THH with K(2) coeff}
The unit map 
	\[\eta \co K(2)\to \THH(\Bt;K(2))\]
	is an equivalence. Consequently, the $P(v_2)$-module $THH_*(\Bt;k(2))$ modulo $v_2$-torsion is freely generated by $1$. 
\end{cor}

\subsection{The $v_2$-Bockstein spectral sequence} \label{BSSv2}
Recall from Section \ref{Bockstein and Adams} that the tower of spectra used to build the Bockstein spectral sequence \eqref{Bockstein and Adams} can be identified as an Adams tower and therefore the Bockstein spectral sequence is multiplicative. 

For $s \geq 4$ recursively define 
\[ \lambda_s := \lambda_{s-3}\mu_3^{p^{s-4}(p-1)}.\] 
For $s \geq 1$ set 
\[ 
r(s,2) = \begin{cases}
          p^s + p^{s-3} + \dots + p^4 + p &  s \equiv 1 \; \text{mod} \; 3 \\
          p^s + p^{s-3} + \dots + p^5 + p^2 & s \equiv 2 \; \text{mod} \; 3 \\
          p^s + p^{s-3} + \dots + p^6 + p^3 & s \equiv 0 \; \text{mod} \; 3 . 
         \end{cases}
\]
\begin{thm}\label{mod p v_1}
Let $\Bt$ be an $E_3$ form of $BP\langle 2 \rangle$. There is an isomorphism of $P(v_2)$-modules 
\[ \THH_*(\Bt; k(2)) \cong P(v_2) \oplus T^2_2,\] 
where 
\begin{align}\label{eq: t22}
T^2_2 \cong \bigoplus_{s \geq 1} P_{r(s,2)}(v_2) \otimes E(\lambda_{s+1}, \lambda_{s+2}) \otimes P(\mu_3^{p^s}) \otimes \mathbb{F}_p\{\lambda_s\mu_3^{jp^{s-1}}: 0 \leq j \leq p-2\}.
\end{align}
\end{thm}
\begin{proof}
We prove by induction on $s \geq 1$ that 
\[ E_{r(s,2)}^{*,*} = P(v_2) \otimes E(\lambda_s, \lambda_{s+1}, \lambda_{s+2}) \otimes P(\mu_3^{p^{s-1}}) \oplus M_s\] 
with 
\[M_s = \bigoplus^{s-1}_{t = 1} P_{r(t,2)}(v_2) \otimes E(\lambda_{t+1}, \lambda_{t+2}) \otimes P(\mu_3^{p^t}) \otimes \mathbb{F}_p\{\lambda_t \mu_3^{jp^{t-1}}: 0 \leq j \leq p-2\},\]
that $\lambda_s$, $\lambda_{s+1}$ and $\lambda_{s+2}$ are infinite cycles,  and that $d_{r(s,2)}(\mu_3^{p^{s-1}}) \dot{=} v_2^{r(s,2)} \lambda_s$.  This implies the statement. 

Since the $v_2^n$ survive to the $E_\infty$-page by Corollary \ref{cor:THH with K(2) coeff}, the classes $\lambda_1$, $\lambda_2$ and $\lambda_3$ are infinite cycles. The class $\mu_3$ needs to support a differential, because otherwise the spectral sequence would collapse at the $E_2$-page by multiplicativity, which is a contradiction to Corollary \ref{cor:THH with K(2) coeff}. 
For bidegree reasons the only possibility is 
\[ d_p(\mu_3) \dot{=} v_2^p \lambda_1.\] 
This proves the base step $s = 1$ of the induction. 
We now assume that the statement holds for some $s \geq 1$. 
We then get 
\[ E_{r(s,2)+1}^{*,*} = P(v_2) \otimes E(\lambda_{s+1}, \lambda_{s+2}, \lambda_s \mu_3^{p^{s-1}(p-1)}) \otimes P(\mu_3^{p^s}) \oplus M_{s+1}.\] 
It now suffices to show that $\lambda_{s+3} =  \lambda_s\mu_3^{p^{s-1}(p-1)}$ is an infinite cycle and that we have a differential 
$d_{r(s+1,2)}(\mu_3^{p^s}) \dot{=} v_2^{r(s+1,2)}\lambda_{s+1}$. 
We cannot have a differential of the form 
\[d_r(\lambda_{s+3}) \dot{=} v_2^n \lambda_{s+1} \lambda_{s+2}\] 
for degree reasons, so $\lambda_{s+3}$ is an infinite cycle. The class $\mu_3^{p^s}$ must support a differential, because otherwise the spectral sequence would collapse at this stage, which is a contradiction to Corollary \ref{cor:THH with K(2) coeff}. For bidegree reasons the only possibility is 
\[ d_{r(s+1,2)}(\mu_3^{p^s}) \dot{=} v_2^{r(s+1,2)} \lambda_{s+1}.\] 
Note that $v_2^{r(s,2)} \lambda_{s+3}$ has the right topological degree, but a too small filtration degree to be the target of a differential on $\mu_3^{p^s}$. This completes the inductive step. 
\end{proof}
\noindent We end with a conjectural answer for $\THH(BP\langle n\rangle;k(m))$ for all $1\le m\le n$. 
\begin{conjecture}\label{conj}
Suppose $1\le m\le n$. Let $\Bn$ be an $E_3$ form of $BP\langle n\rangle$. There is an isomorphism 
\[ 	\THH_*(\Bn ;k(m))\cong  E(\lambda_1,\dots \lambda_{n-m})\otimes  \left ( P(v_m)\oplus T_m^n \right ),\]
where 
\[ T_{m}^{n}=\bigoplus_{s\ge 1} P_{r_n(s,m)}(v_m)\otimes E(\lambda_{n-m+s+1},\dots ,\lambda_{n+s}) \otimes P(\mu_{n+1}^{p^s})\otimes \mathbb{F}_p\{\lambda_{n-m+s}\mu_{n+1}^{p^{\ell p^{s-1}}}:0 \le \ell\le p-2\}\]
and by convention $E(\lambda_1,\dots ,\lambda_{n-m})=\mathbb{F}_p$ when $n=m$. 
The sequence of integers $r_n(s,m)$ is  defined by 
\begin{align*}\label{r(s) for km}
r_n(s,m) = p^{n-m+s}+p^{n-m+s-(m+1)}+\cdots +p^{n+j-m},
\end{align*}
where $j $ is the unique element in $\{1, \dots, m+1\}$ such that $s \equiv j$ mod $m+1$. 

Here the class $\lambda_s$ is defined recursively by the formula 
\[
\lambda_s:=\lambda_{s-(m+1)}\mu_{n+1}^{p^{s-(n+2)}(p-1)} 
\] 
for $s\ge n+2$
and we name the classes in the abutment that are not divisible by $v_n$ by their projection to $\THH_*(\Bn ;H\mathbb{F}_p)$. 
\end{conjecture}

\begin{rem}
When $m=1$ and $n=2$, we observe that this is consistent with Theorem \ref{mod p v_2} where $r_2(s,1)=r(s,1)$. When $m=2$ and $n=2$, we observe that this is consistent with Theorem \ref{mod p v_1} where $r_2(s,2)=r(s,2)$. 
\end{rem}